\documentclass[a4paper,10pt]{article}

\usepackage{amsfonts,amssymb,amsmath,amsthm,xspace}
\usepackage{graphics,pst-tree}
\usepackage{array,multirow}

\newtheorem{theorem}{Theorem}
\newtheorem{proposition}{Proposition}
\newtheorem{definition}{Definition}
\newtheorem{example}{Example}

\newtheorem{remark}{Remark}

\title{Posets and Permutations in the Duplication-Loss Model: Minimal Permutations with $d$ Descents.}
\author{Mathilde Bouvel \thanks{LIAFA, Universit\'e Paris Diderot - Paris 7, Case 7014,
75205 Paris Cedex 13, France, \texttt{mbouvel@liafa.jussieu.fr}} \and Elisa Pergola \thanks{Universit\`a degli Studi di Firenze, Dipartimento di Sistemi e Informatica, v.le Morgagni 65, 50134 Firenze, Italy, \texttt{elisa@dsi.unifi.it}}}

\begin{document}

\maketitle

\begin{abstract}
In this paper, we are interested in the combinatorial analysis of the whole genome duplication - random loss model of genome rearrangement initiated in \cite{CCMR06} and \cite{BR07}. In this model, genomes composed of $n$ genes are modelled by permutations of the set of integers $[1..n]$, that can evolve through duplication-loss steps. It was previously shown that the class of permutations obtained in this model after a given number $p$ of steps is a class of pattern-avoiding permutations of finite basis. The excluded patterns were described as the minimal permutations with $d=2^p$ descents, minimal being intended in the sense of the pattern-involvement relation on permutations. Here, we give a local and simpler characterization of the set $\mathcal{B}_d$ of minimal permutations with $d$ descents. We also provide a more detailed analysis - characterization, bijection and enumeration - of two particular subsets of $\mathcal{B}_d$, namely the patterns in $\mathcal{B}_d$ of size $d+2$ and $2d$.
\end{abstract}

\section{Pattern-avoidance in the duplication-loss model}
\label{section:intro}

The study of genome evolution has been the source of extensive research in computational biology in the last decades. Many models for genome evolution were defined, taking into account various biological phenomema (see \cite{BBCP07}, \cite{CL04}, \cite{journals/tcbb/Labarre06} for recent examples in literature). Among them, the \emph{tandem duplication - random loss model} represents genomes with permutations, that can evolve through \emph{duplication-loss steps} representing the biological phenomenon that duplicates fragments of genomes, and then loses one copy of every duplicated gene. For the original biological motivations, we refer to \cite{CCMR06}. In this first section, we describe the duplication-loss model, and recall some previous results obtained by other authors. We recall some definitions and properties on pattern-avoidance that are necessary to introduce the permutations that will arise from this model and on which we will focus in the rest of the paper.

\subsection{The tandem duplication - random loss model for \\ genome evolution}
\label{section:duplication-loss}

A permutation of size $n$ is a bijective map from $[1..n]$ to itself. We denote by $S_n$ the set of permutations of size $n$. We consider a permutation $\sigma \in S_n$ as the word $\sigma_1 \sigma_2 \ldots \sigma_n$ of $n$ letters on the alphabet $\{1,2,\ldots,n\}$, containing exactly once each letter (we often prefer the word \emph{element} instead of letter). For example, $346251$ represents the permutation $\sigma \in S_6$ such that $\sigma_1 = 3, \sigma_2 = 4, \ldots, \sigma_6 = 1$.

In our model, permutations can be modified by \emph{duplication-loss steps}. Each of these steps is composed of two elementary operations. Firstly, a fragment of consecutive elements of the permutation is duplicated, and the duplicated fragment is inserted immediately after the original copy: this is the \emph{tandem duplication}. After this first operation, any duplicated element appears twice in the sequence of integers (that is no more a permutation at this stage). Then the \emph{random loss} occurs: one copy of every duplicated element is lost, so that we get a permutation at the end of the step. For any duplication-loss step, we call its \emph{width} the number of elements that are duplicated.

\begin{figure}[ht]
\begin{center}
\begin{eqnarray*}
1 \ 2 \ \overbrace{3 \ 4 \ 5 \ 6} \ 7\ & \rightsquigarrow & 1\ 2\ \overbrace{3\ 4\ 5\ 6}\ \overbrace{3 \ 4\ 5\ 6}\ 7  \\
& & \textrm{ (tandem duplication) } \\
 & \rightsquigarrow & 1\ 2\ 3 \ \hspace{-1em} \diagup \, 4\ 5\ 6\ \hspace{-1em} \diagup \, 3 \ 4\ \hspace{-1em} \diagup \, 5\ \hspace{-1em} \diagup \, 6\ 7 \\
& &  \textrm{ (random loss) } \\
& \rightsquigarrow & 1\ 2\ 4\ 5\ 3\ 6\ 7
\end{eqnarray*}
\caption{Example of one step of tandem duplication - random loss of width $4$ \label{ex:duplication-loss}}
\end{center}
\end{figure}

Notice that the duplication-loss model is a particular case of the very general framework for transforming permutations defined in \cite{AAA+04}: the \emph{permuting machines}. A permuting machine takes a permutation in input and performs on it a transformation that satisfies the two properties of independance with respect to the values and of stability with respect to pattern-involvement (see \cite{AAA+04} for more details). These two properties are satisfied by the duplication-loss transformation.

We will consider permutations that are obtained from an identity permutation $12\ldots n$ after a given number $p$ of duplication-loss steps, that is to say that are the output of a combination in series of $p$ permuting machines with input $12\ldots n$. The reason is that these permuations are the ones obtainable at a cost of at most $p$ in the duplication-loss model with a particular \emph{cost function}.

Indeed, various duplication-loss models can be defined depending on the cost function $c \in \overline{\mathbb{R}}^\mathbb{N}$ that is chosen. We will always assume that the cost $c(k)$ of a duplication-loss step is dependant only on the width $k$ of this step. In the original model of Chaudhuri, Chen, Mihaescu and Rao \cite{CCMR06}, the \emph{cost} of a duplication-loss step of width $k$ is $c(k) = \alpha^k$, for a parameter $\alpha \geq 1$. In \cite{BR07}, we consider the cost function defined by $c(k) = 1$ if $k \leq K$, $c(k) = \infty$ otherwise, for a parameter $K \in \overline{\mathbb{N}} \smallsetminus \{0,1\}$. The model we will focus on in what follows has a very simple cost function, namely $c(k) = 1, \forall k$. It is a special case of both the model of \cite{CCMR06} (with $\alpha = 1$) and the model of \cite{BR07} (with $K=\infty$). This particular model is called the \emph{whole genome duplication - random loss model}: indeed, since any step has cost $1$ no matter its width, we can assume \emph{w.l.o.g} that the whole permutation is duplicated at any step.

As said before, we are now going to focus on permutations obtained from an identity permutation $12 \ldots n$ after a certain number $p$ of duplication-loss steps in the whole genome duplication - random loss model, that is to say on permutations obtainable at a cost of at most $p$ in this model. We will describe combinatorial properties of those permutations in Subsection \ref{section:pattern-avoidance}, in terms of pattern-avoidance.

\subsection{Previous results on the duplication-loss model}
\label{section:previous-results}

The permutations obtainable in at most $p$ duplication-loss steps in the whole genome duplication - random loss model were implicitely characterized in \cite{CCMR06}, through Theorem \ref{thm:whole-genome-cost}:
\begin{theorem}
Let $\sigma \in S_n$. In the whole genome duplication - random loss model, $\lceil \log_2 (\textrm{number of maximal increasing substrings of } \sigma ) \rceil$ steps are necessary and sufficient to obtain $\sigma$ from $12\ldots n$.
\label{thm:whole-genome-cost}
\end{theorem}
An \emph{increasing substring} of $\sigma$ is just a sequence of consecutive elements of $\sigma$ that are in increasing order. An increasing substring is maximal if it can be extended neither on the left nor on the right.
\begin{example}
For example, $698413725$ contains $5$ maximal increasing substrings that are $69$, $8$, $4$, $137$ and $25$.
\label{ex:increasing-substring}
\end{example}

In \cite{BR07}, we reformulated Theorem \ref{thm:whole-genome-cost} into Theorem  \ref{thm:whole-genome-nb-desc}, introducing, instead of the number of maximal increasing substrings, the number of descents which is a very well-known statistics on permutations.

\begin{definition}
Given a permutation $\sigma$ of size $n$, we say that there is a \emph{descent} (resp. \emph{ascent}) at position $i$, $1 \leq i \leq n-1$, if $\sigma_i > \sigma_{i+1}$ (resp. $\sigma_i < \sigma_{i+1}$ ). We indicate the number of descents of the permutation $\sigma$ by $desc(\sigma)$. \label{def:descent} 
\end{definition}

\begin{example}
For example, $\sigma = 698413725$ has $4$ descents, namely at positions $2$, $3$, $4$, $7$. \label{ex:descent}
\end{example}

\begin{figure}[ht]
\begin{center}
\psset{unit=0.3cm}
\begin{pspicture}(-4,0)(11,11)
\psgrid[subgriddiv=1,gridwidth=.2pt,griddots=5,gridlabels=0pt](0,0)(9,9)
\rput(0.5,-0.5){{\tiny $1$}}
\rput(1.5,-0.5){{\tiny $2$}}
\rput(2.5,-0.5){{\tiny $3$}}
\rput(3.5,-0.5){{\tiny $4$}}
\rput(4.5,-0.5){{\tiny $5$}}
\rput(5.5,-0.5){{\tiny $6$}}
\rput(6.5,-0.5){{\tiny $7$}}
\rput(7.5,-0.5){{\tiny $8$}}
\rput(8.5,-0.5){{\tiny $9$}}
\rput(-0.5,0.5){{\tiny $1$}}
\rput(-0.5,1.5){{\tiny $2$}}
\rput(-0.5,2.5){{\tiny $3$}}
\rput(-0.5,3.5){{\tiny $4$}}
\rput(-0.5,4.5){{\tiny $5$}}
\rput(-0.5,5.5){{\tiny $6$}}
\rput(-0.5,6.5){{\tiny $7$}}
\rput(-0.5,7.5){{\tiny $8$}}
\rput(-0.5,8.5){{\tiny $9$}}
\pscircle*(0.5,5.5){0.3}
\pscircle*(1.5,8.5){0.3}
\pscircle*(2.5,7.5){0.3}
\pscircle*(3.5,3.5){0.3}
\pscircle*(4.5,0.5){0.3}
\pscircle*(5.5,2.5){0.3}
\pscircle*(6.5,6.5){0.3}
\pscircle*(7.5,1.5){0.3}
\pscircle*(8.5,4.5){0.3}
\psline(0.5,5.5)(1.5,8.5)
\psline(1.5,8.5)(2.5,7.5)
\psline(2.5,7.5)(3.5,3.5)
\psline(3.5,3.5)(4.5,0.5)
\psline(4.5,0.5)(5.5,2.5)
\psline(5.5,2.5)(6.5,6.5)
\psline(6.5,6.5)(7.5,1.5)
\psline(7.5,1.5)(8.5,4.5)
\end{pspicture}
\caption{The grid representation of the permutation $\sigma = 698413725$ \label{fig:grid-representation}}
\end{center}
\end{figure}
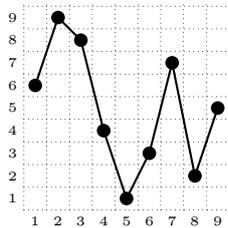

It is often convenient to see permutations through their grid reprensentation defined in \cite{BFP05} and described in Figure \ref{fig:grid-representation}, especially because it gives a better view of descents and ascents.

Obviously, we have:
\begin{remark}
The number of maximal increasing substrings of a permutation $\sigma$ is $desc(\sigma)+1$.
\label{remark:desc}
\end{remark}
More precisely, the positions of the descents and $n$ indicate the positions of the last elements of the maximal increasing substrings of $\sigma$.

These definitions allow us to state Theorem \ref{thm:whole-genome-nb-desc}:
\begin{theorem}
The permutations that can be obtained in at most $p$ steps in the whole genome duplication - random loss model are exactly those whose number of descents is at most $2^p-1$.
\label{thm:whole-genome-nb-desc}
\end{theorem}
\begin{proof}
By Theorem \ref{thm:whole-genome-cost}, the permutations obtainable in at most $p$ steps are exactly those having their number of maximal increasing substrings at most $2^p$, that is to say having at most $2^p-1$ descents by Remark \ref{remark:desc}.
\end{proof}

Generalizing a little, we will focus in the remaining of the paper on the set of permutations with at most $d$ descents, without assuming that $d$ is of the form $d=2^p-1$. We can notice that this corresponds to the set of permutations composed of $d+1$ increasing sequences, separated either by ascents or by descents (a permutation may have more than one such decomposition). In \cite{AMR02} this set is denoted $W(e_1, \ldots, e_{d+1})$ with $\forall i , e_i = +$. In this paper, and as an application of their results, the authors are concerned with properties of $W(e_1, \ldots, e_{d+1})$ in terms of pattern-avoidance, and they prove that this set is a finitely based pattern-avoiding permutation class. Our work can be seen as a more detailled analysis of this particular result.

\subsection{Pattern-avoidance in the duplication-loss model}
\label{section:pattern-avoidance}

We need to recall a few definitions on pattern-avoidance in permutations to proceed.

\begin{definition}
A permutation $\pi \in S_k$ is a \emph{pattern} of a permutation $\sigma \in S_n$ if there is a subsequence of $\sigma$ which is order-isomorphic to $\pi$; i.e., if there is a subsequence $\sigma_{i_1} \sigma_{i_2} \ldots \sigma_{i_k}$ of $\sigma$ (with $1 \leq i_1 < i_2 <\ldots<i_k \leq n$) such that $\sigma_{i_{\ell}} < \sigma_{i_m}$ whenever $\pi_{\ell} < \pi_{m}$. \\
We also say that $\pi$ is \emph{involved} in $\sigma$ and call $\sigma_{i_1} \sigma_{i_2} \ldots \sigma_{i_k}$ an \emph{occurrence} of $\pi$ in $\sigma$. \label{def:pattern}
\end{definition}

\begin{example}
For example $\sigma=1  4  2  5  6  3$ contains the pattern $\pi = 1  3  4  2$; and $1  5  6  3$, $1  4  6  3$, $2 5 6 3$ and $1 4 5 3$ are the occurrences of this pattern in $\sigma$. But $\sigma$ does not contain the pattern $3  2  1$ as no subsequence of size $3$ of $\sigma$  is isomorphic to $3  2  1$, \emph{i.e.}, is decreasing. \label{ex:pattern}
\end{example}

We write $\pi \prec \sigma$ to denote that $\pi$ is a pattern of $\sigma$. We say that a set $\mathcal{C}$ of permutations is \emph{stable for $\prec$} if, for any $\sigma \in \mathcal{C}$, for any $\pi \prec \sigma$, then we also have $\pi \in \mathcal{C}$.

A permutation $\sigma$ that does not contain $\pi$ as a pattern is said to \emph{avoid} $\pi$. The class of all permutations avoiding the patterns $\pi_1, \pi_2 \ldots \pi_k$ is denoted $S(\pi_1, \pi_2, \ldots, \pi_k)$. We say that $S(\pi_1, \pi_2, \ldots, \pi_k)$ is a class of pattern-avoiding permutations of \emph{basis} $\{\pi_1, \pi_2, \ldots, \pi_k\}$. The basis of a class of pattern-avoiding permutations may be finite or infinite. Pattern-avoiding permutation classes considered in the literature (see for example \cite{Bou03}, \cite{Vat05} and their references) are often of finite basis.

Although it may sound a powerful statement, it is simple to understand that:
\begin{proposition}
A set $\mathcal{C}$ of permutations that is stable for $\prec$ is a class of pattern-avoiding permutations. However, its basis might be infinite.
\label{prop:stable-for-prec}
\end{proposition}

\begin{proof}
Consider $\mathcal{C}$ a set of permutations that is stable for $\prec$. Define $\mathcal{B}$ to be the set of minimal permutations that do not belong to $\mathcal{C}$, minimal being intended in the sense of $\prec$. More formally, $\mathcal{B} = \{\sigma \notin \mathcal{C} : \forall \pi \prec \sigma \textrm{ with } \pi \neq \sigma, \pi \in \mathcal{C} \}$. We claim that $\mathcal{C} = S(\mathcal{B})$. Indeed, take $\sigma \notin S(\mathcal{B})$. Then there exists $\pi \in \mathcal{B}$ such that $\pi \prec \sigma$. Since $\pi \in \mathcal{B}$, $\pi \notin \mathcal{C}$ and considering that $\mathcal{C}$ is stable for $\prec$, we deduce from $\pi \prec \sigma$ that $\sigma \notin \mathcal{C}$ either. Conversely, if $\sigma \notin \mathcal{C}$, then either $\sigma \in \mathcal{B}$ (and consequently $\sigma \notin S(\mathcal{B})$) or there exists $ \pi \prec \sigma $ with $\pi \neq \sigma$ such that $\pi \notin \mathcal{C}$. In this second case, by induction we obtain that $\sigma \notin S(\mathcal{B})$.

We conclude that the set $\mathcal{C}$ is a class of pattern-avoiding permutations whose basis $\mathcal{B} = \{\sigma \notin \mathcal{C} : \forall \pi \prec \sigma \textrm{ with } \pi \neq \sigma, \pi \in \mathcal{C} \}$ has no reason \emph{a priori} to be finite.
\end{proof}

In \cite{BR07}, we proved that classes of permutations defined in duplication-loss models, as the permutations obtained in at most a given number $p$ of steps, are classes of pattern-avoiding permutations. We have not always been able to find the basis, even though we have proved in any case we considered that this basis is finite. In this paper, we take into consideration in particular the following result:

\begin{theorem}
The class of permutations obtainable in at most $p$ steps in the whole genome duplication - random loss model is a class of pattern-avoiding permutations whose basis $\mathcal{B}_d$ is finite and is composed of the minimal permutations with $d=2^p$ descents, minimal being intented in the sense of $\prec$.
\label{thm:p-steps}
\end{theorem}

The proof of Theorem \ref{thm:p-steps} we gave in \cite{BR07} is implicite, and for sake of clarity we give below an explicit proof of it.

\begin{proof}
Let us denote by $\mathcal{C}_p$ the class of permutations obtainable in at most $p$ steps in the whole genome duplication - random loss model.

We first prove that $\mathcal{C}_p$ is stable for $\prec$. Consider $\sigma \in \mathcal{C}_p$ of size $n$ and $\pi$ of size $k \leq n$ such that $\pi \prec \sigma$. There is a sequence of at most $p$ duplication-loss steps that transforms $12 \ldots n$ into $\sigma$. By definition, $\pi$ has an occurrence in $\sigma$. In the duplication-loss scenario for $\sigma$, if you keep track only of the elements that form an occurrence of $\pi$, you obtain a sequence of duplication-loss steps moving from $12 \ldots k$ to $\pi$, of no more than $p$ steps. This shows that $\pi \in \mathcal{C}_p$, and consequently that $\mathcal{C}_p$ is stable for $\prec$.

According to Proposition \ref{prop:stable-for-prec}, $\mathcal{C}_p$ is a class of pattern-avoiding permutations whose basis is $ \{\sigma \notin \mathcal{C}_p : \forall \pi \prec \sigma \textrm{ with } \pi \neq \sigma, \pi \in \mathcal{C}_p \}$. Following  Theorem \ref{thm:whole-genome-nb-desc}, we deduce that this basis $\mathcal{B}_d$ of excluded patterns is made of the minimal permutations with $d=2^p$ descents, that is to say the permutations with $2^p$ descents that contain no pattern with $2^p$ descents, except themselves. What is left to prove is that this basis is finite.

It is sufficient to establish an upper bound on the size of the permutations in $\mathcal{B}_d$ to show that $\mathcal{B}_d$ is finite. We postpone this part of the proof to Proposition \ref{prop:d+1<size<2d}, where we show in particular that the permutations of $\mathcal{B}_d$ are of size at most $2d$. A consequence is that the basis $\mathcal{B}_d$ of excluded patterns of $\mathcal{C}_p$ is finite.
\end{proof}

In this paper, we focus on the basis $\mathcal{B}_d$ of excluded patterns appearing in Theorem \ref{thm:p-steps}. More generally, we do not assume that $d$ is a power of $2$ but rather wish to characterize and enumerate the set $\mathcal{B}_d$ of permutations that are the minimal ones in the sense of $\prec$ for the property of having $d$ descents.

\subsection{Outline of the paper}
\label{section:outline}

In this paper, we focus on the sets $\mathcal{B}_d$ of permutations that are the minimal ones in the sense of $\prec$ for the property of having $d$ descents. For the cases $d=2^p$, $\mathcal{B}_d$ is the basis of excluded patterns of the class of permutations obtainable in at most $p$ steps in the whole genome duplication - random loss model.

The work that is presented hereafter is organized as follows. First, we give a local characterization of the permutations of $\mathcal{B}_d$. Indeed, the definition of these permutations as the minimal ones with respect to $\prec$ for the property of having $d$ descents is not very easy to use. We will prove in Section \ref{section:characterization} that the permutations of $\mathcal{B}_d$ are the permutations $\sigma$ whose ascents satisfy a simple and local property: there is an ascent in $\sigma \in S_n$ at position $i$ if and only if $2 \leq i \leq n-2$ and $\sigma_{i-1} \sigma_i \sigma_{i+1} \sigma_{i+2}$ forms an occurrence of either the pattern $2143$ or the pattern $3142$.

This characterization is used to try and count the permutations in $\mathcal{B}_d$. Despite our effort, we did not succeed in this direction, and focused on simpler cases that can be seen as a first step in the enumeration of $\mathcal{B}_d$. First, as explained at the beginning of Section \ref{section:characterization}, Proposition \ref{prop:d+1<size<2d}, the size of the permutations in $\mathcal{B}_d$ is at least $d+1$ and at most $2d$. Obviously there is only one permutation of $\mathcal{B}_d$ of size $d+1$, that is the reverse identity permutation $(d+1)d(d-1) \ldots 321$. For any other size, there is no immediate result. Using a representation of permutations of $\mathcal{B}_d$ as posets (partially ordered sets), we could enumerate the permutations in $\mathcal{B}_d$ having size $2d$ and $d+2$ respectively.

In Section \ref{section:2d}, we prove that the permutations of $\mathcal{B}_d$ having size $2d$ (\emph{i.e.} maximal size) are enumerated by the Catalan numbers: there are $C_d = \frac{1}{d+1} \binom {2d} {d}$ of them. We give two possible proofs of this result. We describe an \emph{``ECO'' generation} (see \cite{BDPP99}) of the permutations of $\mathcal{B}_d$ of size $2d$ whose associated succession rule is known to correspond to the Catalan numbers. More directly, we could provide a simple bijection between Dyck paths of length $2d$ and an adequate representation of the permutations of size $2d$ in $\mathcal{B}_d$.

In Section \ref{section:d+2}, we consider permutations of size $d+2$ (minimal non-trivial case) in $\mathcal{B}_d$. After a combinatorial analysis and some computations, we obtain that there are $s_d = 2^{d+2} - (d+1)(d+2) -2$ such permutations. The sequence $(s_d)$ does not appear in the Online Encyclopedia of Integer Sequences \cite{njas}. However, we realized that the sequence $(\frac{s_d}{2})$ does. This sequence counts the number of non-interval subsets of the set $\{1,2, \ldots , d+1\}$. Section \ref{section:d+2} also gives a bijective proof of the fact that there are twice as many permutations of size $d+2$ in $\mathcal{B}_d$ as non-interval subsets of $\{1,2, \ldots , d+1\}$.

Section \ref{section:conclusion} summarizes some open problems in the study of the sets $\mathcal{B}_d$'s.\\

From here on, by minimal permutation with $d$ descents, we mean a permutation that is minimal in the sense of the pattern-involvement relation $\prec$ for the property of having $d$ descents.

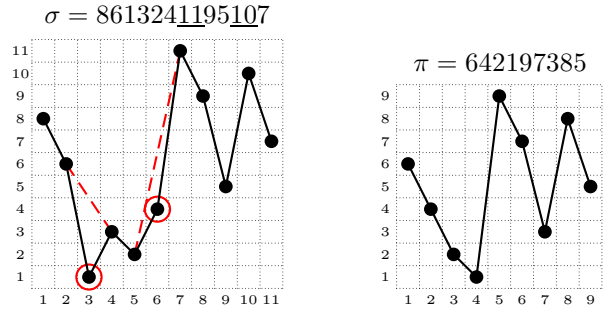
\begin{figure}[ht]
\begin{center}
\psset{unit=0.3cm}
\begin{pspicture}(-5,0)(11,11)
\psgrid[subgriddiv=1,gridwidth=.2pt,griddots=7,gridlabels=0pt](0,0)(11,11)
\rput(0.5,-0.5){{\tiny $1$}}
\rput(1.5,-0.5){{\tiny $2$}}
\rput(2.5,-0.5){{\tiny $3$}}
\rput(3.5,-0.5){{\tiny $4$}}
\rput(4.5,-0.5){{\tiny $5$}}
\rput(5.5,-0.5){{\tiny $6$}}
\rput(6.5,-0.5){{\tiny $7$}}
\rput(7.5,-0.5){{\tiny $8$}}
\rput(8.5,-0.5){{\tiny $9$}}
\rput(9.5,-0.5){{\tiny $10$}}
\rput(10.5,-0.5){{\tiny $11$}}
\rput(-0.5,0.5){{\tiny $1$}}
\rput(-0.5,1.5){{\tiny $2$}}
\rput(-0.5,2.5){{\tiny $3$}}
\rput(-0.5,3.5){{\tiny $4$}}
\rput(-0.5,4.5){{\tiny $5$}}
\rput(-0.5,5.5){{\tiny $6$}}
\rput(-0.5,6.5){{\tiny $7$}}
\rput(-0.5,7.5){{\tiny $8$}}
\rput(-0.5,8.5){{\tiny $9$}}
\rput(-0.5,9.5){{\tiny $10$}}
\rput(-0.5,10.5){{\tiny $11$}}
\psline[linestyle=dashed,linecolor=red](1.5,5.5)(3.5,2.5)
\psline[linestyle=dashed,linecolor=red](4.5,1.5)(6.5,10.5)
\pscircle*(0.5,7.5){0.3}
\pscircle*(1.5,5.5){0.3}
\pscircle*(2.5,0.5){0.3}
\pscircle[linecolor=red](2.5,0.5){0.6}
\pscircle*(3.5,2.5){0.3}
\pscircle*(4.5,1.5){0.3}
\pscircle*(5.5,3.5){0.3}
\pscircle[linecolor=red](5.5,3.5){0.6}
\pscircle*(6.5,10.5){0.3}
\pscircle*(7.5,8.5){0.3}
\pscircle*(8.5,4.5){0.3}
\pscircle*(9.5,9.5){0.3}
\pscircle*(10.5,6.5){0.3}
\psline(0.5,7.5)(1.5,5.5)
\psline(1.5,5.5)(2.5,0.5)
\psline(2.5,0.5)(3.5,2.5)
\psline(3.5,2.5)(4.5,1.5)
\psline(4.5,1.5)(5.5,3.5)
\psline(5.5,3.5)(6.5,10.5)
\psline(6.5,10.5)(7.5,8.5)
\psline(7.5,8.5)(8.5,4.5)
\psline(8.5,4.5)(9.5,9.5)
\psline(9.5,9.5)(10.5,6.5)
\rput(5.5,12){$\sigma = 8 6 1 3 2 4 \underline{11} 9 5 \underline{10} 7$}
\end{pspicture}\begin{pspicture}(-5,0)(11,13)
\psgrid[subgriddiv=1,gridwidth=.2pt,griddots=7,gridlabels=0pt](0,0)(9,9)
\rput(0.5,-0.5){{\tiny $1$}}
\rput(1.5,-0.5){{\tiny $2$}}
\rput(2.5,-0.5){{\tiny $3$}}
\rput(3.5,-0.5){{\tiny $4$}}
\rput(4.5,-0.5){{\tiny $5$}}
\rput(5.5,-0.5){{\tiny $6$}}
\rput(6.5,-0.5){{\tiny $7$}}
\rput(7.5,-0.5){{\tiny $8$}}
\rput(8.5,-0.5){{\tiny $9$}}
\rput(-0.5,0.5){{\tiny $1$}}
\rput(-0.5,1.5){{\tiny $2$}}
\rput(-0.5,2.5){{\tiny $3$}}
\rput(-0.5,3.5){{\tiny $4$}}
\rput(-0.5,4.5){{\tiny $5$}}
\rput(-0.5,5.5){{\tiny $6$}}
\rput(-0.5,6.5){{\tiny $7$}}
\rput(-0.5,7.5){{\tiny $8$}}
\rput(-0.5,8.5){{\tiny $9$}}
\pscircle*(0.5,5.5){0.3}
\pscircle*(1.5,3.5){0.3}
\pscircle*(2.5,1.5){0.3}
\pscircle*(3.5,0.5){0.3}
\pscircle*(4.5,8.5){0.3}
\pscircle*(5.5,6.5){0.3}
\pscircle*(6.5,2.5){0.3}
\pscircle*(7.5,7.5){0.3}
\pscircle*(8.5,4.5){0.3}
\psline(0.5,5.5)(1.5,3.5)
\psline(1.5,3.5)(2.5,1.5)
\psline(2.5,1.5)(3.5,0.5)
\psline(3.5,0.5)(4.5,8.5)
\psline(4.5,8.5)(5.5,6.5)
\psline(5.5,6.5)(6.5,2.5)
\psline(6.5,2.5)(7.5,7.5)
\psline(7.5,7.5)(8.5,4.5)
\rput(4.5,10){$\pi = 6 4 2 1 9 7 3 8 5$}
\end{pspicture}
\caption{Permutations $\sigma$ and $\pi$ of Example \ref{ex:d-desc}. \label{fig:d-desc} }
\end{center}
\end{figure}

Example \ref{ex:d-desc}, which is illustrated on Figure \ref{fig:d-desc}, should clarify the notion of minimal permutation with $d$ descents.

\begin{example}
Permutation $\sigma = 8 6 1 3 2 4 \underline{11} 9 5 \underline{10} 7$ has $6$ descents but is not minimal with $6$ descents. Indeed, the elements $1$ and $4$ (that are circled on Figure \ref{fig:d-desc}) can be removed from $\sigma$ without changing the number of descents. 

Doing this, we obtain permutation $\pi = 6 4 2 1 9 7 3 8 5$ which is minimal with $6$ descents: it is impossible to remove an element from it while preserving the number of descents equal to $6$. 

However, $\pi$ is not of minimal \emph{size} among the permutation with $6$ descents: $\pi$ has size $9$ whereas permutation $7654321$ has $6$ descents but size $7$.
\label{ex:d-desc}
\end{example}

\section{A characterization for minimal permutations with $d$ descents}
\label{section:characterization}

The aim of this section is to provide a more practical characterization of minimal permutations with $d$ descents, by finding necessary and sufficient conditions on permutations for being minimal with $d$ descents. First, we provide a necessary condition on the size of those permutations with Proposition \ref{prop:d+1<size<2d}.
\begin{proposition}
Let $\sigma$ be a minimal permutation with $d$ descents. Then every ascent of $\sigma$ is immediately preceded and immediately followed by a descent, and the size $n$ of $\sigma$ satisfies $d+1 \leq n \leq 2d$.
\label{prop:d+1<size<2d}
\end{proposition}

\begin{proof}
Consider a permutation $\sigma \in \mathcal{B}_d$, and denote by $n$ the size of $\sigma$. By minimality in the sense of $\prec$, $\sigma$ has exactly $d$ descents. To create a permutation with $d$ descents, you need at least $d+1$ elements, and with $d+1$ elements, the only permutation with $d$ descents you can create is $(d+1) d (d-1) \ldots 21$, which is minimal. Therefore, $n \geq d+1$.

It is also easily seen that $\sigma$ does neither start nor end with an ascent, otherwise the permutation obtained by removing the first or the last element of $\sigma$ would have the same number $d$ of descents, contradicting that $\sigma$ is minimal with $d$ descents. In the same way, $\sigma$ cannot have two consecutive ascents $\sigma_{i-1} \sigma_i$ and $\sigma_i \sigma_{i+1}$, otherwise we would get the same contradiction removing $\sigma_i$, since this removal does not change the number of descents.

This proves that a minimal permutation with $d$ descents is composed of non-empty sequences of descents, separated by isolated ascents. A longest possible permutation with $d$ descents so obtained has $d$ isolated descents, separated by $d-1$ isolated ascents, and consequently has $2d$ elements. We then get that the size of $\sigma$ is at most $2d$: $n \leq 2d$.
\end{proof}

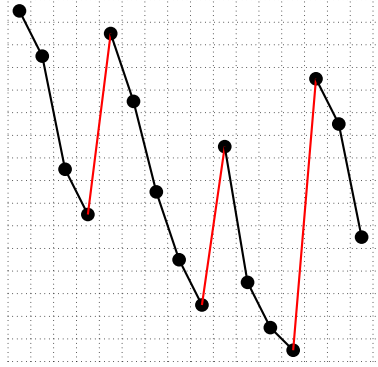
\begin{figure}[ht]
\begin{center}
\psset{unit=0.3cm}
\begin{pspicture}(-4,0)(18,18)
\psgrid[subgriddiv=1,gridwidth=.2pt,griddots=5,gridlabels=0pt](0,0)(16,16)
\pscircle*(0.5,15.5){0.3}
\pscircle*(1.5,13.5){0.3}
\pscircle*(2.5,8.5){0.3}
\pscircle*(3.5,6.5){0.3}
\pscircle*(4.5,14.5){0.3}
\pscircle*(5.5,11.5){0.3}
\pscircle*(6.5,7.5){0.3}
\pscircle*(7.5,4.5){0.3}
\pscircle*(8.5,2.5){0.3}
\pscircle*(9.5,9.5){0.3}
\pscircle*(10.5,3.5){0.3}
\pscircle*(11.5,1.5){0.3}
\pscircle*(12.5,0.5){0.3}
\pscircle*(13.5,12.5){0.3}
\pscircle*(14.5,10.5){0.3}
\pscircle*(15.5,5.5){0.3}
\psline(0.5,15.5)(1.5,13.5)
\psline(1.5,13.5)(2.5,8.5)
\psline(2.5,8.5)(3.5,6.5)
\psline[linecolor=red](3.5,6.5)(4.5,14.5)
\psline(4.5,14.5)(5.5,11.5)
\psline(5.5,11.5)(6.5,7.5)
\psline(6.5,7.5)(7.5,4.5)
\psline(7.5,4.5)(8.5,2.5)
\psline[linecolor=red](8.5,2.5)(9.5,9.5)
\psline(9.5,9.5)(10.5,3.5)
\psline(10.5,3.5)(11.5,1.5)
\psline(11.5,1.5)(12.5,0.5)
\psline[linecolor=red](12.5,0.5)(13.5,12.5)
\psline(13.5,12.5)(14.5,10.5)
\psline(14.5,10.5)(15.5,5.5)
\end{pspicture}
\caption{Decomposition of a minimal permutation with $d$ descents into non-empty sequences of descents separated by isolated ascents \label{fig:decomposition}}
\end{center}
\end{figure}

The decomposition of a minimal permutation with $d$ descents into non-empty sequences of descents separated by isolated ascents that is described in the proof of Proposition  \ref{prop:d+1<size<2d} is illustrated in Figure \ref{fig:decomposition}. This decomposition can be carried further to give a necessary and sufficient condition on permutations for being minimal with $d$ descents. This characterization is described in Theorem \ref{thm:characterization}.

\begin{theorem}
A permutation $\sigma$ is minimal with $d$ descents if and only if it has exactly $d$ descents and its ascents $\sigma_i \sigma_{i+1}$ are such that $2 \leq i \leq n-2$ and $\sigma_{i-1} \sigma_i \sigma_{i+1} \sigma_{i+2}$ forms an occurrence of either the pattern $2143$ or the pattern $3142$.
\label{thm:characterization}
\end{theorem}

\begin{proof}
Let $\sigma$ be a minimal permutation with $d$ descents. In the decomposition of $\sigma$ into non-empty sequences of descents separated by isolated ascents -- illustrated in Figure \ref{fig:decomposition} -- it appears clearly that an ascent $\sigma_i \sigma_{i+1}$ is necessarily such that $2 \leq i \leq n-2$, with $\sigma_{i-1} \sigma_i$ and $\sigma_{i+1} \sigma_{i+2}$ being descents.

Now, consider an ascent $\sigma_i \sigma_{i+1}$. The previous remarks lead to $\sigma_{i-1} > \sigma_i$, $\sigma_i < \sigma_{i+1}$ and $\sigma_{i+1} > \sigma_{i+2}$.

Let us assume that $\sigma_{i-1} > \sigma_{i+1}$. Then the permutation obtained from $\sigma$ by the removal of $\sigma_i$ has as many descents as $\sigma$ (and one ascent less), contradicting the minimality of $\sigma$. Consequently, $\sigma_{i-1} < \sigma_{i+1}$. Similarly, if $\sigma_{i} > \sigma_{i+2}$, the removal of $\sigma_{i+1}$ from $\sigma$ does not change the number of descents, contradicting the minimality of $\sigma$. So $\sigma_{i} < \sigma_{i+2}$ (see Figure \ref{fig:forbidden}).

\begin{figure}[ht]
\begin{center}

\psset{unit=0.3cm}
\begin{pspicture}(-4,0)(22,6)
\pscircle*(-2.5,3.5){0.3}
\pscircle*(-1.5,0.5){0.3}
\pscircle*(-0.5,2.5){0.3}
\pscircle*(0.5,1.5){0.3}
\psline(-2.5,3.5)(-1.5,0.5)
\psline(-1.5,0.5)(-0.5,2.5)
\psline(-0.5,2.5)(0.5,1.5)

\pscircle[linecolor=red](-1.5,0.5){0.6}
\psline[linestyle=dotted,linecolor=red](-2.5,3.5)(-0.5,2.5)

\pscircle*(2.5,2.5){0.3}
\pscircle*(3.5,1.5){0.3}
\pscircle*(4.5,3.5){0.3}
\pscircle*(5.5,0.5){0.3}
\psline(2.5,2.5)(3.5,1.5)
\psline(3.5,1.5)(4.5,3.5)
\psline(4.5,3.5)(5.5,0.5)

\pscircle[linecolor=red](4.5,3.5){0.6}
\psline[linestyle=dotted,linecolor=red](3.5,1.5)(5.5,0.5)

\rput(-8.3,3){Removing} \pscircle[linecolor=red](-5,3){0.4}
\rput(-7.85,2){leads to}
\rput(-8.85,1){the descent} \psline[linestyle=dotted,linecolor=red](-6.8,1.4)(-5,0.4)

\rput(2.5,5){Forbidden configurations}

\psline[linestyle=dashed](10.5,-0.5)(10.5,6.5)

\pscircle*(12.5,2.5){0.3}
\pscircle*(13.5,0.5){0.3}
\pscircle*(14.5,3.5){0.3}

\pscircle*(15.5,1.5){0.3}
\psline(12.5,2.5)(13.5,0.5)
\psline(13.5,0.5)(14.5,3.5)
\psline(14.5,3.5)(15.5,1.5)

\pscircle*(17.5,1.5){0.3}
\pscircle*(18.5,0.5){0.3}
\pscircle*(19.5,3.5){0.3}
\pscircle*(20.5,2.5){0.3}
\psline(17.5,1.5)(18.5,0.5)
\psline(18.5,0.5)(19.5,3.5)
\psline(19.5,3.5)(20.5,2.5)

\rput(19.5,5){The only possible configurations}

\end{pspicture}
\caption{The elements $\sigma_{i-1} \sigma_i \sigma_{i+1} \sigma_{i+2}$ around an ascent $\sigma_i \sigma_{i+1}$ in a permutation $\sigma$ which is minimal with $d$ descents \label{fig:forbidden}}
\end{center}
\end{figure}
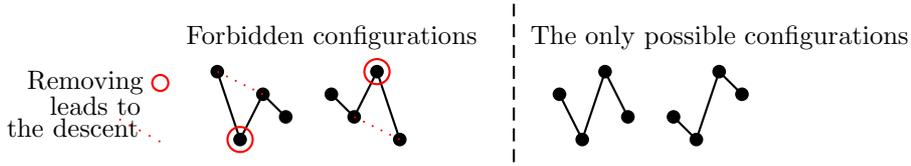

At this point, we have the five following inequalities: $\sigma_{i-1} > \sigma_i$, $\sigma_i < \sigma_{i+1}$, $\sigma_{i+1} > \sigma_{i+2}$, $\sigma_{i-1} < \sigma_{i+1}$ and $\sigma_{i} < \sigma_{i+2}$. Thanks to them it is possible to check that the sequence $\sigma_{i-1} \sigma_i \sigma_{i+1} \sigma_{i+2}$ is an occurrence of either the pattern $2143$ or the pattern $3142$.

Conversely, consider a permutation $\sigma$ with $d$ descents whose ascents $\sigma_i \sigma_{i+1}$ are such that $2 \leq i \leq n-2$ and $\sigma_{i-1} \sigma_i \sigma_{i+1} \sigma_{i+2}$ forms an occurrence of either the pattern $2143$ or the pattern $3142$. This implies that $\sigma$ has the shape of non-empty sequences of descents separated by isolated ascents. And it is a simple matter to prove that the removal of any element of $\sigma$ makes the number of descents decrease by one -- there are three cases to consider: the removed element may be either the first element of an ascent, or the second element of an ascent, or it may be between two descents. This proves that $\sigma$ is a minimal permutation with $d$ descents.
\end{proof}

We thought this characterization could help us to enumerate the minimal permutations with $d$ descents. Although we did not reach this goal, we still obtain partial results when we studied minimal permutation with $d$ descents and of a given size $n$. For $n = d+1$, we already proved that there is only one such permutation. For $n= d+2$ and $n=2d$, the next two sections describe the enumeration we obtained. In both cases, we will use a \emph{partially ordered set} (or \emph{poset}) representation of permutations, that comes directly from the characterization of minimal permutations with $d$ descents in Theorem \ref{thm:characterization}.

\paragraph{Representation of minimal permutations with $d$ descents with posets}

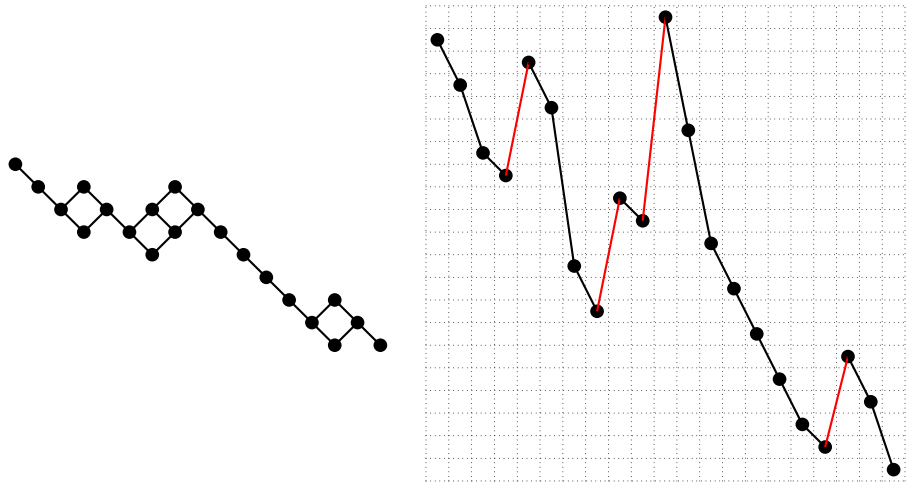
\begin{figure}[ht]
\begin{center}

\psset{unit=0.30cm}
\begin{pspicture}(0,-6)(40,16)
\pscircle*(1,9){0.3}
\psline(1,9)(2,8)
\pscircle*(2,8){0.3}
\psline(2,8)(3,7)
\pscircle*(3,7){0.3}
\psline(3,7)(4,6)
\psline(3,7)(4,8)
\pscircle*(4,6){0.3}
\psline(4,6)(5,7)
\pscircle*(4,8){0.3}
\psline(4,8)(5,7)
\pscircle*(5,7){0.3}
\psline(5,7)(6,6)
\pscircle*(6,6){0.3}
\psline(6,6)(7,5)
\psline(6,6)(7,7)
\pscircle*(7,5){0.3}
\psline(7,5)(8,6)
\pscircle*(7,7){0.3}
\psline(7,7)(8,6)
\psline(7,7)(8,8)
\pscircle*(8,6){0.3}
\psline(8,6)(9,7)
\pscircle*(8,8){0.3}
\psline(8,8)(9,7)
\pscircle*(9,7){0.3}
\psline(9,7)(10,6)
\pscircle*(10,6){0.3}
\psline(10,6)(11,5)
\pscircle*(11,5){0.3}
\psline(11,5)(12,4)
\pscircle*(12,4){0.3}
\psline(12,4)(13,3)
\pscircle*(13,3){0.3}
\psline(13,3)(14,2)
\pscircle*(14,2){0.3}
\psline(14,2)(15,1)
\psline(14,2)(15,3)
\pscircle*(15,1){0.3}
\psline(15,1)(16,2)
\pscircle*(15,3){0.3}
\psline(15,3)(16,2)
\pscircle*(16,2){0.3}
\psline(16,2)(17,1)
\pscircle*(17,1){0.3}

\psgrid[subgriddiv=1,gridwidth=.2pt,griddots=5,gridlabels=0pt](19,-5)(40,16)
\pscircle*(19.5,14.5){0.3}
\pscircle*(20.5,12.5){0.3}
\pscircle*(21.5,9.5){0.3}
\pscircle*(22.5,8.5){0.3}
\pscircle*(23.5,13.5){0.3}
\pscircle*(24.5,11.5){0.3}
\pscircle*(25.5,4.5){0.3}
\pscircle*(26.5,2.5){0.3}
\pscircle*(27.5,7.5){0.3}
\pscircle*(28.5,6.5){0.3}
\pscircle*(29.5,15.5){0.3}
\pscircle*(30.5,10.5){0.3}
\pscircle*(31.5,5.5){0.3}
\pscircle*(32.5,3.5){0.3}
\pscircle*(33.5,1.5){0.3}
\pscircle*(34.5,-0.5){0.3}
\pscircle*(35.5,-2.5){0.3}
\pscircle*(36.5,-3.5){0.3}
\pscircle*(37.5,0.5){0.3}
\pscircle*(38.5,-1.5){0.3}
\pscircle*(39.5,-4.5){0.3}
\psline(19.5,14.5)(20.5,12.5)
\psline(20.5,12.5)(21.5,9.5)
\psline(21.5,9.5)(22.5,8.5)
\psline[linecolor=red](22.5,8.5)(23.5,13.5)
\psline(23.5,13.5)(24.5,11.5)
\psline(24.5,11.5)(25.5,4.5)
\psline(25.5,4.5)(26.5,2.5)
\psline[linecolor=red](26.5,2.5)(27.5,7.5)
\psline(27.5,7.5)(28.5,6.5)
\psline[linecolor=red](28.5,6.5)(29.5,15.5)
\psline(29.5,15.5)(30.5,10.5)
\psline(30.5,10.5)(31.5,5.5)
\psline(31.5,5.5)(32.5,3.5)
\psline(32.5,3.5)(33.5,1.5)
\psline(33.5,1.5)(34.5,-0.5)
\psline(34.5,-0.5)(35.5,-2.5)
\psline(35.5,-2.5)(36.5,-3.5)
\psline[linecolor=red](36.5,-3.5)(37.5,0.5)
\psline(37.5,0.5)(38.5,-1.5)
\psline(38.5,-1.5)(39.5,-4.5)

\end{pspicture}

\caption{A poset representing a set of minimal permutations with $16$ descents and $4$ ascents (consequently of size $21$) containing, among others the permutation $20 \ 18 \ 15 \ 14 \ 19 \ 17 \ 10 \ 8 \ 13 \ 12 \ 21 \ 16 \ 11 \ 9 \ 7 \ 5 \ 3 \ 2 \ 6 \ 4 \ 1 $ whose grid representation is also given \label{fig:poset}}
\end{center}
\end{figure}

Consider a set of all the permutations of a given size $n$, that are minimal with $d$ descents, and having their descents and ascents in the same positions. In all these permutations, the elements are locally ordered in the same way, even around the ascents, because of Theorem \ref{thm:characterization}. We can give a representation of this whole set of permutations by a \emph{partially ordered set} (or \emph{poset}) indicating the necessary conditions on the relative order of the elements between them. For a descent, we just have a link from the first and greatest element to the second and smallest one. For any ascent $\sigma_i \sigma_{i+1}$, the elements $\sigma_{i-1} \sigma_i \sigma_{i+1} \sigma_{i+2}$ form a diamond-shaped structure with $\sigma_{i+1}$ on the top, $\sigma_i$ on the bottom, $\sigma_{i-1}$ on the left and $\sigma_{i+2}$ on the right. See Figure \ref{fig:poset} for an example. By Theorem \ref{thm:characterization}, any labelling of the elements of the poset respecting its ordering constraints is a minimal permutation with $d$ descents.

We will say that a permutation $\sigma$ satisfies the diamond property when each of its ascent $\sigma_i \sigma_{i+1}$ is such that $\sigma_{i-1} \sigma_i \sigma_{i+1} \sigma_{i+2} $ forms a diamond, that is to say is an occurrence of either $2143$ or $3142$.

\section{Enumeration of minimal permutations with $d$ descents and of size $2d$}
\label{section:2d}

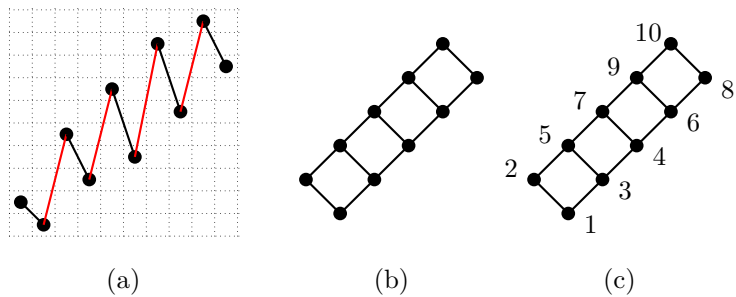
\begin{figure}[ht]
\begin{center}
\psset{unit=0.3cm}
\begin{pspicture}(0,-2)(33,10)
\psgrid[subgriddiv=1,gridwidth=.2pt,griddots=5,gridlabels=0pt](0,0)(10,10)
\pscircle*(0.5,1.5){0.3}
\pscircle*(1.5,0.5){0.3}
\pscircle*(2.5,4.5){0.3}
\pscircle*(3.5,2.5){0.3}
\pscircle*(4.5,6.5){0.3}
\pscircle*(5.5,3.5){0.3}
\pscircle*(6.5,8.5){0.3}
\pscircle*(7.5,5.5){0.3}
\pscircle*(8.5,9.5){0.3}
\pscircle*(9.5,7.5){0.3}
\psline(0.5,1.5)(1.5,0.5)
\psline[linecolor=red](1.5,0.5)(2.5,4.5)
\psline(2.5,4.5)(3.5,2.5)
\psline[linecolor=red](3.5,2.5)(4.5,6.5)
\psline(4.5,6.5)(5.5,3.5)
\psline[linecolor=red](5.5,3.5)(6.5,8.5)
\psline(6.5,8.5)(7.5,5.5)
\psline[linecolor=red](7.5,5.5)(8.5,9.5)
\psline(8.5,9.5)(9.5,7.5)
\rput(5,-2){(a)}

\pscircle*(14.5,1){0.3}
\pscircle*(13,2.5){0.3}
\pscircle*(16,2.5){0.3}
\pscircle*(14.5,4){0.3}
\pscircle*(17.5,4){0.3}
\pscircle*(16,5.5){0.3}
\pscircle*(19,5.5){0.3}
\pscircle*(17.5,7){0.3}
\pscircle*(20.5,7){0.3}
\pscircle*(19,8.5){0.3}
\psline(13,2.5)(14.5,1)
\psline(14.5,4)(16,2.5)
\psline(16,5.5)(17.5,4)
\psline(17.5,7)(19,5.5)
\psline(19,8.5)(20.5,7)
\psline(13,2.5)(14.5,4)
\psline(14.5,4)(16,5.5)
\psline(16,5.5)(17.5,7)
\psline(17.5,7)(19,8.5)
\psline(14.5,1)(16,2.5)
\psline(16,2.5)(17.5,4)
\psline(17.5,4)(19,5.5)
\psline(19,5.5)(20.5,7)
\rput(16.75,-2){(b)}

\pscircle*(24.5,1){0.3}
\rput(25.5,0.5){$1$}
\pscircle*(23,2.5){0.3}
\rput(22,3){$2$}
\pscircle*(26,2.5){0.3}
\rput(27,2){$3$}
\pscircle*(24.5,4){0.3}
\rput(23.5,4.5){$5$}
\pscircle*(27.5,4){0.3}
\rput(28.5,3.5){$4$}
\pscircle*(26,5.5){0.3}
\rput(25,6){$7$}
\pscircle*(29,5.5){0.3}
\rput(30,5){$6$}
\pscircle*(27.5,7){0.3}
\rput(26.5,7.5){$9$}
\pscircle*(30.5,7){0.3}
\rput(31.5,6.5){$8$}
\pscircle*(29,8.5){0.3}
\rput(28,9){$10$}
\psline(23,2.5)(24.5,1)
\psline(24.5,4)(26,2.5)
\psline(26,5.5)(27.5,4)
\psline(27.5,7)(29,5.5)
\psline(29,8.5)(30.5,7)
\psline(23,2.5)(24.5,4)
\psline(24.5,4)(26,5.5)
\psline(26,5.5)(27.5,7)
\psline(27.5,7)(29,8.5)
\psline(24.5,1)(26,2.5)
\psline(26,2.5)(27.5,4)
\psline(27.5,4)(29,5.5)
\psline(29,5.5)(30.5,7)
\rput(26.75,-2){(c)}
\end{pspicture}
\caption{(a) A minimal permutation $\sigma =2 \ 1 \ 5 \ 3 \ 7 \ 4 \ 9 \ 6 \ 10 \ 8 $ with $d=5$ descents and of size $2d =10$, (b) the poset representing the set of all minimal permutations with $d=5$ descents and of size $2d =10$ and (c) the authorized labelling of the subsequent poset associated with $\sigma$ \label{fig:2d-ex}}
\end{center}
\end{figure}

The minimal permutations with $d$ descents that have size $2d$ are, because of minimality, of a very particular shape. Indeed, they cannot have two consecutive ascents as usual, but neither can they have two consecutive descents, otherwise it would be impossible to reach size $2d$. Consequently, they all result from of an alternation of isolated descents and isolated ascents, of course starting and ending with a descent. An example is given in Figure \ref{fig:2d-ex}(a).

A consequence is that all minimal permutations of size $2d$ with $d$ descents have their descents and ascents in the same position, so that a unique poset represents the set of all minimal permutations with $d$ descents having size $2d$. This poset has the shape of a ladder with $d$ steps: it is a sequence of $d-1$ diamonds, two consecutive diamonds being linked by an edge. These diamonds correspond to the ascents in the permutations, that are separated by one descent only in this case. See Figure \ref{fig:2d-ex}(b) for an example.

The paragraph on poset representation at the end of Section \ref{section:characterization} justifies Proposition \ref{prop:2d-characterization}:
\begin{proposition}
The minimal permutations with $d$ descents and of size $2d$ correspond exactly to the labellings of the ladder poset with $d$ steps with the integers $\{1,2,\ldots,2d\}$ that respect its ordering constraints. \label{prop:2d-characterization}
\end{proposition}

An example of this correspondance is given in Figure \ref{fig:2d-ex}(c).

The poset representation allows to see at once some properties of minimal permutations with $d$ descents having size $2d$. For example, such a permutation always has $1$ as its second element and $2d$ as its next to last element.

The main result of this section is :

\begin{theorem}
The minimal permutations with $d$ descents and of size $2d$ are enumerated by the Catalan numbers $C_d = \frac{1}{d+1} \binom {2d} {d}$. \label{thm:2d}
\end{theorem}

\paragraph{Proof of Theorem \ref{thm:2d} by an analytical method}
A possible way to prove Theorem \ref{thm:2d} is to use the ECO method, presented in details in \cite{BDPP99}. In our case, the idea developed by this method is to build all the authorized labellings of the ladder poset with $d$ steps from all the authorized labellings of the ladder poset with $d-1$ steps without creating twice the same labelling.

In its original form, the ECO method builds combinatorial objects of size $d$ from those of size $d-1$, through a process of \emph{local expansion}, whereby the objects are modified only by the addition of an elementary block of object. In our case, in order to get a labelling of size $d$, the \emph{local expansions} might modify many labels in the labelling of size $d-1$, but the relative order of these labels between them will remain the same. In this sense, we can consider that the expansion is still local.

In the ECO method, the combinatorial objects (labellings of the ladder poset with $d$ steps in our case) receive \emph{labels}. The label of an object is the number of its \emph{children}, that is to say the number of objects that are obtained from it in the local expansion process. Those children can again receive a label by the same method. The infinite tree in which any permutation is the father of its children is called the \emph{generating tree} of the combinatorial class.

With the ECO labelling of the combinatorial objects, we derive a \emph{succession rule} or \emph{rewriting rule} that describes the production (in terms of labels) of the possible labels of these objects, together with a starting point. There is a simple succession rule that is associated with some combinatorial classes enumerated by the Catalan numbers (for example with Dyck paths \cite{BDPP99}):
$$
\left\{
\begin{array}{l}
  (2) \\
  (k) \rightsquigarrow (2) (3) \cdots (k) (k+1)\\
\end{array}
\right.
$$

A possible way of proving that authorized labellings of the ladder poset with $d$ steps are enumerated by the Catalan numbers is to find an ECO construction for this class whose associated succession rule is the one above.

The ECO labels that are given to authorized labellings of the ladder posets with $d$ steps for this purpose are $(2d-\sigma_{2d} +1)$, $\sigma_{2d}$ being the label of the rightmost element of the poset. Notice also that $2d$ is the label of the uppermost element of the poset, and that this element is also the second rightmost one. 

Consider an authorized labelling $\sigma$ of the ladder poset with $d$ steps that has ECO label $(k)$. Its children are the labellings of the ladder poset with $d+1$ steps obtained by adding a new step on the right, this new step of the ladder being labelled with $2d+2$ for the top element, and $i$ for the rightmost one, for $2d+2-k \leq i \leq 2d+1$. The elements $j$ in $\sigma$ with $j \geq i$ are turned into to $j+1$ to maintain both the relative order of the elements of $\sigma$ and the property that these new labellings use all the integers of $[1..2d+2]$ exactly once.

Since $k = 2d-\sigma_{2d} +1$, it is easy to check that all the labellings obtained in this way are authorized, and that all of them are obtained. We can now focus on the ECO labels of the children (of size $d+1$) of an authorized labelling of size $d$ with ECO label $(k)$. There are of course $k$ of them whose ECO labels are, by the above formula, $(2(d+1) - i +1)$ with $2d+2-k \leq i \leq 2d+1$, that is to say the children of a labelling with ECO label $(k)$ have labels $(2), (3), \ldots (k), (k+1)$.

The starting point for this ECO construction is the ladder poset with one step provided with its only authorized labelling $21$, and whose ECO label is $(2)$.

To sum up, the succession rule obtained for this ECO construction of authorized labelling of the ladder posets is 
$$
\left\{
\begin{array}{l}
  (2) \\
  (k) \rightsquigarrow (2) (3) \cdots (k) (k+1)\\
\end{array}
\right.
$$
and this succession rule corresponds to combinatorial classes enumerated by Catalan numbers.

\begin{figure}[ht]
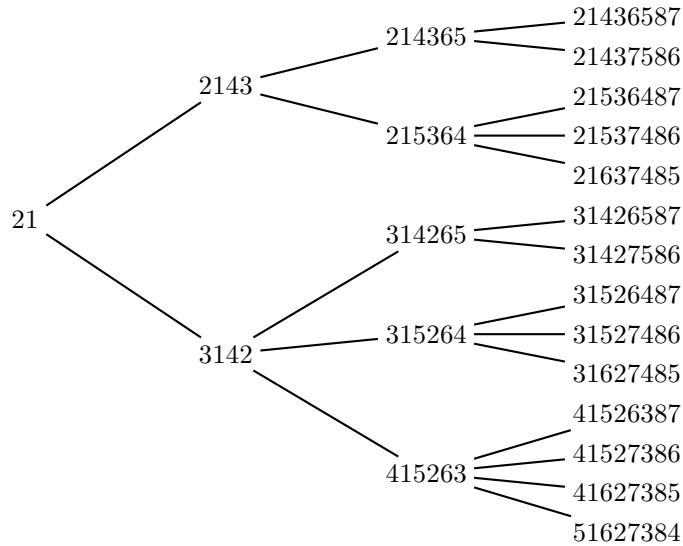

\begin{center}
\psset{levelsep=75pt,treesep=0.3cm,nodesep=3pt}
\pstree[treemode=R]{\Tr{$21$}}{
\pstree{\Tr{$2143$}}{
\pstree{\Tr{$214365$}}{
\Tr{$21436587$} \Tr{$21437586$}}
\pstree{\Tr{$215364$}}{
\Tr{$21536487$} \Tr{$21537486$} \Tr{$21637485$}}
}
\pstree{\Tr{$3142$}}{
\pstree{\Tr{$314265$}}{
\Tr{$31426587$} \Tr{$31427586$}}
\pstree{\Tr{$315264$}}{
\Tr{$31526487$} \Tr{$31527486$} \Tr{$31627485$}}
\pstree{\Tr{$415263$}}{
\Tr{$41526387$} \Tr{$41527386$} \Tr{$41627385$} \Tr{$51627384$}}
}
}

\caption{The first four levels of the generating tree associated with the ECO construction of authorized labellings of the ladder posets described above \label{fig:generating-tree}}
\end{center}
\end{figure}

Figure \ref{fig:generating-tree} shows the beginning of the generating tree associated with this ECO construction. To improve the understanding of this tree, we do not represent labellings of ladder posets in its nodes, but rather the minimal permutations with $d$ descents of size $2d$ associated with them.

\paragraph{Proof of Theorem \ref{thm:2d} by bijection}
It is well known that Dyck paths of length $2d$ are enumerated by the Catalan numbers $C_d = \frac{1}{d+1} \binom {2d} {d}$. Let us recall the definition of Dyck paths.
\begin{definition}
A Dyck path of length $2d$ is a path in $\mathbb{N} \times \mathbb{N}$ starting at $(0,0)$ and ending at $(2d,0)$, with steps going up (of coordinate $(1,1)$) and steps going down (of coordinate $(1,-1)$).
\end{definition}
As it is a path in $\mathbb{N} \times \mathbb{N}$, a Dyck path never goes under the $x$-axis. We can also notice that a Dyck path has as many steps going up as those going down, and that any prefix of a Dyck path contains at least as many steps going up as those going down. This is actually a characterization of Dyck paths.

\begin{figure}[ht]
\begin{center}
\psset{unit=0.3cm}
\begin{pspicture}(0,1)(22,8)
\psgrid[subgriddiv=1,gridwidth=.2pt,griddots=3,gridlabels=0pt](0,2)(10,5)
\pscircle(0,2){0.2}
\rput(0.5,1){{\footnotesize $1$}}
\psline(0,2)(1,3)
\pscircle(1,3){0.2}
\rput(1.5,1){{\footnotesize $2$}}
\psline(1,3)(2,4)
\pscircle(2,4){0.2}
\rput(2.5,1){{\footnotesize $3$}}
\psline(2,4)(3,3)
\pscircle(3,3){0.2}
\rput(3.5,1){{\footnotesize $4$}}
\psline(3,3)(4,4)
\pscircle(4,4){0.2}
\rput(4.5,1){{\footnotesize $5$}}
\psline(4,4)(5,5)
\pscircle(5,5){0.2}
\rput(5.5,1){{\footnotesize $6$}}
\psline(5,5)(6,4)
\pscircle(6,4){0.2}
\rput(6.5,1){{\footnotesize $7$}}
\psline(6,4)(7,3)
\pscircle(7,3){0.2}
\rput(7.5,1){{\footnotesize $8$}}
\psline(7,3)(8,2)
\pscircle(8,2){0.2}
\rput(8.5,1){{\footnotesize $9$}}
\psline(8,2)(9,3)
\pscircle(9,3){0.2}
\rput(9.5,1){{\footnotesize $10$}}
\psline(9,3)(10,2)
\pscircle(10,2){0.2}

\pscircle*(14.5,1){0.3}
\rput(15.5,0.5){$1$}
\pscircle*(13,2.5){0.3}
\rput(12,3){$3$}
\pscircle*(16,2.5){0.3}
\rput(17,2){$2$}
\pscircle*(14.5,4){0.3}
\rput(13.5,4.5){$6$}
\pscircle*(17.5,4){0.3}
\rput(18.5,3.5){$4$}
\pscircle*(16,5.5){0.3}
\rput(15,6){$7$}
\pscircle*(19,5.5){0.3}
\rput(20,5){$5$}
\pscircle*(17.5,7){0.3}
\rput(16.5,7.5){$8$}
\pscircle*(20.5,7){0.3}
\rput(21.5,6.5){$9$}
\pscircle*(19,8.5){0.3}
\rput(18,9){$10$}
\psline(13,2.5)(14.5,1)
\psline(14.5,4)(16,2.5)
\psline(16,5.5)(17.5,4)
\psline(17.5,7)(19,5.5)
\psline(19,8.5)(20.5,7)
\psline(13,2.5)(14.5,4)
\psline(14.5,4)(16,5.5)
\psline(16,5.5)(17.5,7)
\psline(17.5,7)(19,8.5)
\psline(14.5,1)(16,2.5)
\psline(16,2.5)(17.5,4)
\psline(17.5,4)(19,5.5)
\psline(19,5.5)(20.5,7)
\end{pspicture}
\caption{An example of the bijection between minimal permutations with $d$ descents of size $2d$ (seen as authorized labellings of the ladder poset with $d$ steps) and Dyck paths with $2d$ steps \label{fig:bijection-catalan}}
\end{center}
\end{figure}
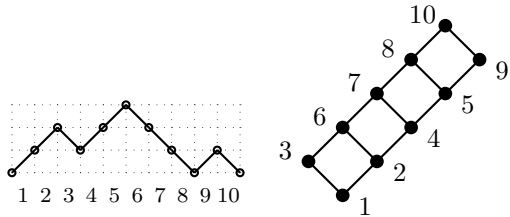

We provide a bijection between Dyck paths of length $2d$ and authorized labellings of the ladder poset with $d$ steps with the integers $\{1,2,\ldots,2d\}$. The bijection is simple. Starting from a Dyck path $\mathcal{D}$ of length $2d$, we number its steps with the integers from $1$ to $2d$, from left to right. Then, we label the lower line of the ladder with the numbers of the steps of $\mathcal{D}$ going up and its upper line with the numbers of the steps of $\mathcal{D}$ going down. An example is shown in Figure \ref{fig:bijection-catalan}.

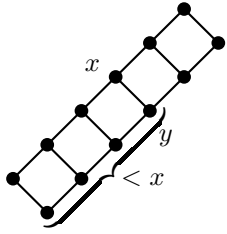
\begin{figure}[ht]
\begin{center}
\psset{unit=0.3cm}
\begin{pspicture}(26,0)(40,7.5)
\pscircle*(28.5,0){0.3}
\pscircle*(27,1.5){0.3}
\pscircle*(30,1.5){0.3}
\pscircle*(28.5,3){0.3}
\pscircle*(31.5,3){0.3}
\pscircle*(30,4.5){0.3}
\pscircle*(33,4.5){0.3}
\rput(33.7,3.3){$y$}
\pscircle*(31.5,6){0.3}
\rput(30.5,6.5){$x$}
\pscircle*(34.5,6){0.3}
\pscircle*(33,7.5){0.3}
\pscircle*(36,7.5){0.3}
\pscircle*(34.5,9){0.3}
\psline(28.5,0)(30,1.5)
\psline(27,1.5)(28.5,0)
\psline(27,1.5)(28.5,3)
\psline(28.5,3)(30,1.5)
\psline(30,4.5)(31.5,3)
\psline(31.5,6)(33,4.5)
\psline(33,7.5)(34.5,6)
\psline(34.5,9)(36,7.5)
\psline(28.5,3)(30,4.5)
\psline(30,4.5)(31.5,6)
\psline(31.5,6)(33,7.5)
\psline(33,7.5)(34.5,9)
\psline(30,1.5)(31.5,3)
\psline(31.5,3)(33,4.5)
\psline(33,4.5)(34.5,6)
\psline(34.5,6)(36,7.5)

\rput(32.7,1.5){$< x$}
\rput{45}(31.2,1.8){$\underbrace{~ \ \qquad \qquad \quad \ ~}$}
\end{pspicture}
\caption{A condition for a labelling of the ladder poset to be authorized \label{fig:authorized-labelling}}
\end{center}
\end{figure}

The application we described is actually a bijection between Dyck paths and the authorized labellings of the ladder posets, corresponding to the permutations we are interested in. The reason is simple. It is sufficient to notice that a labelling of the ladder poset with $d$ steps is authorized if and only if any $i$-th element $x$ on the upper line has at least $i$ smaller elements on the lower line (the element $y$ on the lower line that is linked to $x$ by a step on the ladder, and all the elements below $y$). See Figure \ref{fig:authorized-labelling} for a better understanding of this statement. In the same way, a path with $d$ steps going up and $d$ step going down is a Dyck path if and only if any $i$-th step going down has at least $i$ steps going up before it.

\section{Enumeration of minimal permutations with $d$ descents and of size $d+2$}
\label{section:d+2}

In Section \ref{section:2d}, we enumerated the minimal permutations with $d$ descents and of size $2d$, that is to say of \emph{maximal} possible size. We have already proved that the \emph{minimal} possible size for a minimal permutation with $d$ descents is $d+1$ and shown that there is only one such permutation, namely the reversed identity $(d+1)d \ldots 21$. In this section, we will focus on the minimal permutations with $d$ descents and of size $d+2$, \emph{i.e.} the minimal non-trivial case, and give a closed formula for their enumeration through Theorem \ref{thm:d+2}.

\begin{theorem}
The minimal permutations with $d$ descents and of size $d+2$ are enumerated by the sequence $(s_d)$ defined as follows: $s_d = 2^{d+2} - (d+1)(d+2) -2$. \label{thm:d+2}
\end{theorem}

We provide two possible proofs for Theorem \ref{thm:d+2}. Both of them are based on the poset representation of minimal permutations with $d$ descents and of size $d+2$, that consequently have a unique ascent. The first one is straightforward with this decomposition, but implies rather complex computations. The second proof is more complicated but it does not involve such technicalities: it consists in a correspondance between non-interval subsets of $\{1, 2, \ldots, d+1\}$ and minimal permutations with $d$ descents of size $d+2$, each non-interval subset being associated with exactly two distinct permutations.

\paragraph{Proof of Theorem \ref{thm:d+2} by a computational method} Let us recall that a minimal permutation $\sigma$ with $d$ descents and of size $d+2$ has a unique ascent, between two sequences of descents, and that the elements surrounding the ascent are organized in a diamond in the poset representation of the permutation.

Let us denote by $i$ and $k$ the elements of the ascent, $i<k$, by $j$ the element preceeding $i$ in $\sigma$, and by $h$ the element following $k$. In the permutation $\sigma$, the subsequence $jikh$ forms an occurence of either the pattern $2143$ (if $j<h$) or the pattern $3142$ (if $j>h$). This defines two types of minimal permutations with $d$ descents of size $d+2$. We denote by $N_1$ the number of those permutations for which $j<h$ and by $N_2$ the number of those having $j>h$.

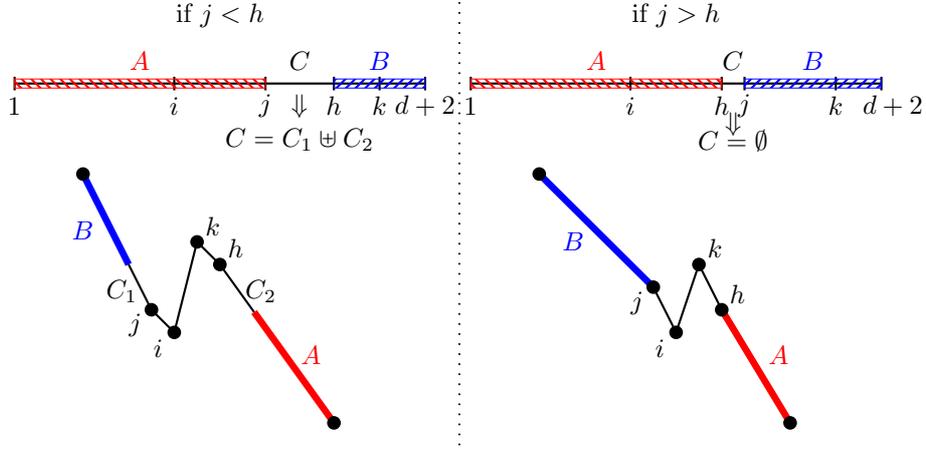
\begin{figure}[ht]
\begin{center}
\psset{unit=0.3cm}
\begin{pspicture}(0,0)(40,20)

\rput(9,19){if $j<h$}

\psline(0,16)(18,16)
\psline(0,15.7)(0,16.3)
\rput(0,15){$1$}
\psline(18,15.7)(18,16.3)
\rput(18,15){$d+2$}

\psline(7,15.7)(7,16.3)
\rput(7,15){$i$}
\psline(11,15.7)(11,16.3)
\rput(11,15){$j$}
\psline(14,15.7)(14,16.3)
\rput(14,15){$h$}
\psline(16,15.7)(16,16.3)
\rput(16,15){$k$}

\psframe[linecolor=red,fillstyle=vlines,hatchcolor=red,linewidth=0,hatchangle=45,hatchsep=0.2](0,15.8)(11,16.2)
\rput(5.5,17){${\red A}$}
\psframe[linecolor=blue,fillstyle=vlines,hatchcolor=blue,linewidth=0,hatchangle=-45,hatchsep=0.2](14,15.8)(18,16.2)
\rput(16,17){${\blue B}$}
\rput(12.5,17){$C$}

\rput(12.5,13.5){$ C = C_1 \uplus C_2$}
\rput(12.5,15){$\Downarrow$}

\psline[linecolor=blue,linewidth=0.3](3,12)(5,8)
\psline[linecolor=red,linewidth=0.3](10.5,5.9)(14,1)
\pscircle*(3,12){0.3}
\rput(3,9.5){${\blue B}$}
\psline(5,8)(6,6)
\rput(4.7,6.8){$C_1$}
\pscircle*(6,6){0.3}
\rput(5.3,5.3){$j$}
\psline(6,6)(7,5)
\pscircle*(7,5){0.3}
\rput(6.3,4.3){$i$}
\psline(7,5)(8,9)
\pscircle*(8,9){0.3}
\rput(8.7,9.7){$k$}
\psline(8,9)(9,8)
\pscircle*(9,8){0.3}
\rput(9.7,8.7){$h$}
\psline(9,8)(10.5,5.9)
\rput(10.8,6.8){$C_2$}
\rput(13,4){${\red A}$}

\pscircle*(14,1){0.3}

\psline[linestyle=dotted](19.5,0)(19.5,20)

\rput(29,19){if $j>h$}

\psline(20,16)(38,16)
\psline(20,15.7)(20,16.3)
\rput(20,15){$1$}
\psline(38,15.7)(38,16.3)
\rput(38.5,15){$d+2$}

\psline(27,15.7)(27,16.3)
\rput(27,15){$i$}
\psline(31,15.7)(31,16.3)
\rput(31,15){$h$}
\psline(32,15.7)(32,16.3)
\rput(32,15){$j$}
\psline(36,15.7)(36,16.3)
\rput(36,15){$k$}

\psframe[linecolor=red,fillstyle=vlines,hatchcolor=red,linewidth=0,hatchangle=45,hatchsep=0.2](20,15.8)(31,16.2)
\rput(25.5,17){${\red A}$}
\psframe[linecolor=blue,fillstyle=vlines,hatchcolor=blue,linewidth=0,hatchangle=-45,hatchsep=0.2](32,15.8)(38,16.2)
\rput(35,17){${\blue B}$}
\rput(31.5,17){$C$}

\rput(31.5,13.5){$ C = \emptyset$}
\rput(31.5,14.2){$\Downarrow$}

\psline[linecolor=blue,linewidth=0.3](23,12)(28,7)
\psline[linecolor=red,linewidth=0.3](31,6)(34,1)
\pscircle*(23,12){0.3}
\rput(24.5,9){${\blue B}$}
\pscircle*(28,7){0.3}
\rput(27.3,6.3){$j$}
\psline(28,7)(29,5)
\pscircle*(29,5){0.3}
\rput(28.3,4.3){$i$}
\psline(29,5)(30,8)
\pscircle*(30,8){0.3}
\rput(30.7,8.7){$k$}
\psline(30,8)(31,6)
\pscircle*(31,6){0.3}
\rput(31.7,6.7){$h$}
\rput(33.5,4){${\red A}$}
\pscircle*(34,1){0.3}

\end{pspicture}
\caption{The two types of minimal permutations with $d$ descents and of size $d+2$, with the decomposition used for their enumeration \label{fig:d+2}}
\end{center}
\end{figure}

We first compute $N_1$. In order to characterize a minimal permutation $\sigma$ with $d$ descents, of size $d+2$, and having its diamond of the type $2143$, you first need to establish the values of $j$, $i$, $k$ and $h$ satisfying the constraints $1 \leq i < j < h < k \leq d+2$. Then (see left part of Figure \ref{fig:d+2}), the elements greater than $h$ (except $k$) are necessarily placed before $j$, in decreasing order, forming the sequence of descents $B$. Similarily, the elements smaller than $j$ (except $i$) have to come after $h$ in $\sigma$, again in decreasing order, to form the sequence of descents $A$. The set $C$ of elements between $j$ and $h$ must be partitioned into two parts $C_1$ and $C_2$, possibly empty, the elements of $C_1$ being placed in decreasing order between $B$ and $j$, those of $C_2$ between $h$ and $A$. There are $2^{card(C)} = 2^{h-j-1}$ such partitions of $C$ into $C_1 \uplus C_2$.

To sum up, a minimal permutation with $d$ descents, of size $d+2$, and having its diamond of the type $2143$ is determined by the values of its $i$, $j$, $h$, and $k$, with $1 \leq i < j < h < k \leq d+2$, and a partition of the set $C$ of elements between $j$ and $h$ into $C_1 \uplus C_2$. This characterization allows us to compute $N_1$:
\begin{eqnarray*}
N_1 & = & \sum_{i=1}^{d-1} \sum_{j=i+1}^{d} \sum_{h=j+1}^{d+1} \sum_{k=h+1}^{d+2} 2^{h-j-1} \\
& = & \sum_{i=1}^{d-1} \sum_{j=i+1}^{d} \sum_{h=j+1}^{d+1} (d+2-h) 2^{h-j-1}\\
& = & \sum_{i=1}^{d-1} \sum_{j=i+1}^{d} \sum_{m=0}^{d-j} (d+1-m-j) 2^m \\
& = & \sum_{i=1}^{d-1} \sum_{j=i+1}^{d} \Big[ \sum_{m=0}^{d-j} (d+1-j) 2^m - \sum_{m=0}^{d-j} m 2^m \Big]\\
& = & \sum_{i=1}^{d-1} \sum_{j=i+1}^{d} \Big[ (d+1-j)(2^{d-j+1}-1) - 2^{d-j+1} (d-j-1) -2 \Big] \\
& = & \sum_{i=1}^{d-1} \sum_{j=i+1}^{d} 2^{d-j+2} - (d-j+2) -1\\
& = & \sum_{i=1}^{d-1} \sum_{n=2}^{d-i+1} 2^{n} - n -1\\
& = & \sum_{i=1}^{d-1} 2^{d-i+2} - \frac{(d-i+1)(d-i+2)}{2} - (d-i) -3 \\
& = & \sum_{p=3}^{d+1} 2^p - \frac{p(p-1)}{2} - (p-2) -3 \\
& = & \sum_{p=3}^{d+1} 2^p -\frac{1}{2} p^2 -\frac{1}{2} p -1 \\
& = & 2^{d+2} -\frac{1}{2} \frac{(d+1)(d+2)(2d+3)}{6} -\frac{1}{2} \frac{(d+1)(d+2)}{2} -d -3 \\
& = & 2^{d+2} - \frac{(d+1)(d+2)(d+3)}{6} -d-3
\end{eqnarray*}

For the minimal permutations $\sigma$ with $d$ descents and of size $d+2$, whose diamond is of type $3142$, the analysis is simpler (this case is illustrated on the right side of Figure \ref{fig:d+2}). Indeed, following the previous notations, to characterize such a permutation, you must again choose $i$, $j$, $h$ and $k$ with the constraint that $1 \leq i < h < j < k$, but not every such choice is acceptable. Namely, consider the set of elements between $h$ and $j$. Those elements cannot be before $j$ in $\sigma$, since they are smaller than $j$. But neither can they go after $h$ since they are greater than $h$. Consequently, there cannot be any element between $h$ and $j$, and $h=j-1$. Now, once $i$, $j$ and $k$ are established, the permutation $\sigma$ is completly characterized. The elements greater than $j$ (except $k$) necessarily form a sequence $B$ of descents before $j$, and those smaller than $j-1$ (except $i$) form a sequence $A$ of descents after $h = j-1$. The computation of $N_2$ is then straigthforward:

\begin{eqnarray*}
N_2 & = & \sum_{i=1}^{d-1} \sum_{j=i+2}^{d+1} \sum_{k=j+1}^{d+2} 1  \\
& = & \sum_{i=1}^{d-1} \sum_{j=i+2}^{d+1} d+2-j \\
& = & \sum_{i=1}^{d-1} \sum_{m=1}^{d-i} m \\
& = & \sum_{i=1}^{d-1} \frac{(d-i)(d-i+1)}{2}\\
& = & \sum_{n=1}^{d-1} \frac{n(n+1)}{2}
\end{eqnarray*}
\begin{eqnarray*}
& = & \frac{1}{2} \Big[ \frac{d(d-1)(2d-1)}{6} + \frac{d(d+1)}{2}\Big] \\
& = & \frac{d(d-1)(d+1)}{6}
\end{eqnarray*}

The total number of minimal permutations with $d$ descents of size $d+2$ is now simply obtained by the final computation:
\begin{eqnarray*}
N = N_1 + N_2 & = & 2^{d+2} - \frac{(d+1)(d+2)(d+3)}{6} -d-3 + \frac{d(d-1)(d+1)}{6}\\
& = & 2^{d+2} - (d+1)^2 -d-3\\
& = & 2^{d+2} -(d+1)(d+2) -2
\end{eqnarray*}

This achieves the computational proof of Theorem \ref{thm:d+2}. We now turn to a bijective proof of it.

\paragraph{Proof of Theorem \ref{thm:d+2} by bijection}
A \emph{non-interval subset} of $\{1,2, \ldots, d+1\}$ is a non-empty subset of $\{1,2, \ldots , d+1\}$ that is not an interval. For example, the non-interval subsets of $\{1, \ldots, 4\}$ are: $\{1,3\}$, $\{1,4\}$, $\{2,4\}$, $\{1,2,4\}$ and $\{1,3,4\}$. Non-interval subsets of $\{1,2, \ldots, d+1\}$ are easy to enumerate, as shown in Proposition \ref{prop:non-interval}.

\begin{proposition}
The number of non-interval subsets of the set $\{1,2, \ldots, d+1\}$ is $2^{d+1} - \frac{(d+1)(d+2)}{2} -1$.
\label{prop:non-interval}
\end{proposition}

\begin{proof}
There are $2^{d+1}$ subsets of $\{1,2, \ldots, d+1\}$, one being the empty set. So we only need to prove that there are $\frac{(d+1)(d+2)}{2}$ subsets of $\{1,2, \ldots, d+1\}$ that \emph{are} (non-empty) intervals. It is simple to see that there are $i$ interval subsets of $\{1,2, \ldots, d+1\}$ whose greatest element is $i$, namely the intervals $[j..i]$ for $1 \leq j \leq i$. And since $ \sum_{i=1}^{d+1} i = \frac{(d+1)(d+2)}{2} \textrm{ ,}$ the proof of Proposition \ref{prop:non-interval} is completed.
\end{proof}

Notice that the sequence $(2^{d+1} - \frac{(d+1)(d+2)}{2} -1)_d$ is registered in the Online Encyclopedia of Integer Sequences \cite{njas} as [A002662]. To prove Theorem \ref{thm:d+2}, we need to show that there are twice as many minimal permutations with $d$ descents and of size $d+2$ as non-interval subsets of $\{1,2, \ldots , d+1\}$. For this purpose, we partition the set of minimal permutations with $d$ descents and of size $d+2$ into two subsets $S^1$ and $S^2$, and show bijections between $S^1$ (resp. $S^2$) and the set of non-interval subsets of $\{1,2, \ldots , d+1\}$, denoted $\mathcal{NI}$.

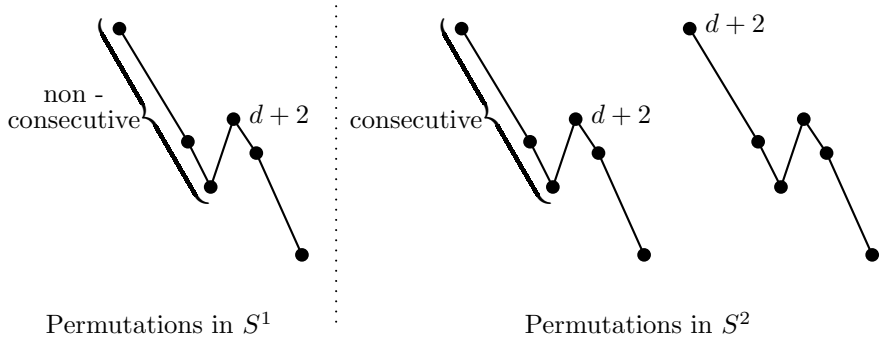
\begin{figure}[ht]
\begin{center}
\psset{unit=0.3cm}
\begin{pspicture}(0,-2)(40,12)

\pscircle*(6,11){0.3}
\psline(6,11)(9,6)
\pscircle*(9,6){0.3}
\psline(9,6)(10,4)
\pscircle*(10,4){0.3}
\psline(10,4)(11,7)
\pscircle*(11,7){0.3}
\rput(13,7.2){$d+2$}
\psline(11,7)(12,5.5)
\pscircle*(12,5.5){0.3}
\psline(12,5.5)(14,1)
\pscircle*(14,1){0.3}

\rput{-60}(7.3,7.2){$\underbrace{~ \ ~ \ \qquad \qquad \qquad \ ~}$}
\rput(4,8){non -}
\rput(4,7){consecutive}

\rput(7.75,-2){Permutations in $S^1$}

\psline[linestyle=dotted](15.5,-2)(15.5,12)

\pscircle*(21,11){0.3}
\psline(21,11)(24,6)
\pscircle*(24,6){0.3}
\psline(24,6)(25,4)
\pscircle*(25,4){0.3}
\psline(25,4)(26,7)
\pscircle*(26,7){0.3}
\rput(28,7.2){$d+2$}
\psline(26,7)(27,5.5)
\pscircle*(27,5.5){0.3}
\psline(27,5.5)(29,1)
\pscircle*(29,1){0.3}

\rput{-60}(22.3,7.2){$\underbrace{~ \ ~ \ \qquad \qquad \qquad \ ~}$}
\rput(19,7){consecutive}

\pscircle*(31,11){0.3}
\rput(33,11.2){$d+2$}
\psline(31,11)(34,6)
\pscircle*(34,6){0.3}
\psline(34,6)(35,4)
\pscircle*(35,4){0.3}
\psline(35,4)(36,7)
\pscircle*(36,7){0.3}
\psline(36,7)(37,5.5)
\pscircle*(37,5.5){0.3}
\psline(37,5.5)(39,1)
\pscircle*(39,1){0.3}

\rput(28.75,-2){Permutations in $S^2$}

\end{pspicture}
\caption{The shapes of the permutations in the sets $S^1$ and $S^2$ \label{fig:S1-and-S2}}
\end{center}
\end{figure}

The set $S^1$ contains the minimal permutations $\sigma$ with $d$ descents and of size $d+2$ such that (1) $d+2$ is the element at the top of the ascent of $\sigma$, and (2) the first sequence of descents of $\sigma$ is not composed of elements that are consecutive. The set $S^2$ contains all the other minimal permutations with $d$ descents and of size $d+2$. Figure \ref{fig:S1-and-S2} shows the shapes of the permutations in $S^1$ and in $S^2$.

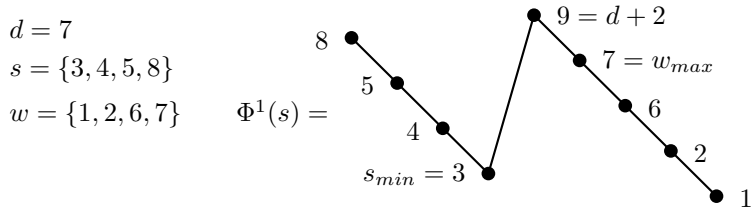
\begin{figure}[ht]
\begin{center}
\psset{unit=0.3cm}
\begin{pspicture}(0,0)(35,10)

\rput[bl](1,8){$d=7$}
\rput[bl](1,6){$s = \{3,4,5,8\}$}
\rput[bl](1,4){$w= \{1,2,6,7\}$}

\rput[br](15,4){$\Phi^1(s) =$}

\pscircle*(16,8){0.3}
\rput[br](15,7.5){$8$}
\psline(16,8)(18,6)
\pscircle*(18,6){0.3}
\rput[br](17,5.5){$5$}
\psline(18,6)(20,4)
\pscircle*(20,4){0.3}
\rput[br](19,3.5){$4$}
\psline(20,4)(22,2)
\pscircle*(22,2){0.3}
\rput[br](21,1.5){$s_{min} = 3$}
\psline(22,2)(24,9)
\pscircle*(24,9){0.3}
\rput[bl](25,8.5){$9 = d+2$}
\psline(24,9)(26,7)
\pscircle*(26,7){0.3}
\rput[bl](27,6.5){$7 = w_{max}$}
\psline(26,7)(28,5)
\pscircle*(28,5){0.3}
\rput[bl](29,4.5){$6$}
\psline(28,5)(30,3)
\pscircle*(30,3){0.3}
\rput[bl](31,2.5){$2$}
\psline(30,3)(32,1)
\pscircle*(32,1){0.3}
\rput[bl](33,0.5){$1$}

\end{pspicture}
\caption{Definition of the bijection $\Phi^1$ on an example \label{fig:phi1-ex}}
\end{center}
\end{figure}

We first describe the simple bijection $\Phi^1$ between $\mathcal{NI}$ and $S^1$. Consider a non-interval subset $s$ of $\{1,2, \ldots, d+1\}$. Let us denote by $w$ the set of ``wholes'' associated with $s$: $w= \{1,2, \ldots , d+1\} \setminus s$. Now we set $\Phi^1(s)$ to be the permutation consisting of the elements of $s$ in decreasing order, followed by $d+2$ and then by the elements of $w$ in decreasing order. This definition is illustrated in Figure \ref{fig:phi1-ex}.

\begin{proposition}
The application $\Phi^1$ defines a bijection between $\mathcal{NI}$ and $S^1$. \label{prop:phi1}
\end{proposition}

\begin{proof}
Let $s$ be a non-interval subset of $\{1,2, \ldots, d+1\}$, and let $w$ be the associated set of wholes $w= \{1,2, \ldots , d+1\} \setminus s$.

We start by proving that $\Phi^1(s) \in S^1$. Since $s \in \mathcal{NI}$, $s$ contains at least two elements, and $w$ at least one. Consequently,  $\Phi^1(s)$ consists of two non-empty sequences of descents separated by one ascent, and we just need to check the diamond property around its ascent to prove that $\Phi^1(s)$ is a minimal permutation with $d$ descents and of size $d+2$. In our case, proving this diamond property is the same as showing that the smallest element $s_{min}$ of $s$ is smaller than the bigger element $w_{max}$ of $w$. Since $s$ is not an interval, there is at least one element of $w$ that is bigger than $s_{min}$, and consequently $s_{min} < w_{max}$. Finally, considering again that $s$ is not an interval, we get that $\Phi^1(s) \in S^1$.

Now -- given that among the minimal permutations with $d$ descents and of size $d+2$, the permutations of $S^1$ are defined as those whose elements in the first sequence of descents do not form an interval -- it should now be clear that $\Phi^1$ is a bijection between $\mathcal{NI}$ and $S^1$.
\end{proof}

\begin{figure}[ht]
\begin{center}
\psset{unit=0.3cm}
\begin{pspicture}(0,-22)(40,12)

\rput(6,9){Type $A$}

\pscircle*(7,6){0.3}
\psline(7,6)(8,4)
\pscircle*(8,4){0.3}
\psline(8,4)(9,7)
\pscircle*(9,7){0.3}
\rput(11,7.2){$d+2$}
\psline(9,7)(10,5.5)
\pscircle*(10,5.5){0.3}
\rput(12,5.7){$d+1$}
\psline(10,5.5)(12,1)
\pscircle*(12,1){0.3}

\rput{-60}(7,4.5){$\underbrace{~  \qquad \ ~}$}
\rput(3,4){consecutive}


\rput(27,10.5){Type $B$}

\pscircle*(21,11){0.3}
\rput(22,11.2){$d$}
\psline(21,11)(23.8,7)
\pscircle*(23.8,7){0.3}
\psline(23.8,7)(25,5)
\pscircle*(25,5){0.3}
\psline(25,5)(26,8)
\pscircle*(26,8){0.3}
\rput(28,8.2){$d+2$}
\psline(26,8)(27,6.5)
\pscircle*(27,6.5){0.3}
\rput(29,6.7){$d+1$}
\psline(27,6.5)(29,1)
\pscircle*(27.5,5){0.3}
\pscircle*(29,1){0.3}
\rput(28,0.8){$1$}

\rput{-55}(22.5,7.5){$\underbrace{~ \ \qquad \qquad \qquad \ ~}$}
\rput(19,7){consecutive}

\rput{110}(29,3.25){$\underbrace{~  \qquad \qquad \ ~}$}
\rput(32.5,4){consecutive}
\rput(32.5,3){or empty}


\rput(13.5,-2.5){Type $C$}

\pscircle*(10,-3){0.3}
\psline(10,-3)(14,-7)
\pscircle*(14,-7){0.3}
\psline(14,-7)(16,-9)
\pscircle*(16,-9){0.3}
\psline(16,-9)(17,-4)
\pscircle*(17,-4){0.3}
\rput(19,-3.5){$d+2$}
\psline(17,-4)(19,-5)
\pscircle*(19,-5){0.3}
\rput(21,-4.5){$d+1$}
\psline(19,-5)(21,-6)
\pscircle*(21,-6){0.3}
\rput(20.4,-6.7){$d$}
\psline(21,-6)(25,-8)
\pscircle*(25,-8){0.3}
\psline(25,-8)(27,-9)
\pscircle*(27,-9){0.3}
\psline(27,-9)(31,-11)
\pscircle*(31,-11){0.3}
\rput(30.5,-11.7){$1$}

\rput{-45}(12.5,-6.5){$\underbrace{~ \ \quad \qquad \qquad \qquad \ ~}$}
\rput(9,-7){consecutive}

\rput{153}(23.3,-6.4){$\underbrace{~ \ ~ \ \quad \qquad \ ~}$}
\rput[bl](23,-6){consecutive and \textbf{not} empty}

\rput{153}(29.3,-9.4){$\underbrace{~ \ ~ \ \quad \qquad \ ~}$}
\rput[bl](29,-9){consecutive or empty}

\rput(11,-14.5){Type $D$}

\pscircle*(5,-13){0.3}
\rput[bl](5.8,-13.2){$d+2$}
\psline(5,-13)(6,-15)
\pscircle*(6,-15){0.3}
\psline(6,-15)(8,-19)
\pscircle*(8,-19){0.3}
\psline(8,-19)(9,-21)
\pscircle*(9,-21){0.3}
\psline(9,-21)(11,-17)
\pscircle*(11,-17){0.3}
\psline(11,-17)(12,-18)
\pscircle*(12,-18){0.3}
\psline(12,-18)(14,-20)
\pscircle*(14,-20){0.3}

\rput{315}(13,-18){$\overbrace{~ \ ~ \ \quad \qquad \ ~}$}
\rput[bl](13.5,-18){consecutive}

\rput(28,-14.5){Type $E$}

\pscircle*(22,-13){0.3}
\rput[bl](22.8,-13.2){$d+2$}
\psline(22,-13)(23,-15)
\pscircle*(23,-15){0.3}
\psline(23,-15)(25,-19)
\pscircle*(25,-19){0.3}
\psline(25,-19)(26,-21)
\pscircle*(26,-21){0.3}
\psline(26,-21)(28,-17)
\pscircle*(28,-17){0.3}
\psline(28,-17)(29,-18)
\pscircle*(29,-18){0.3}
\psline(29,-18)(31,-20)
\pscircle*(31,-20){0.3}

\rput{315}(30,-18){$\overbrace{~ \ ~ \ \quad \qquad \ ~}$}
\rput[bl](30.5,-18){non-consecutive}

\end{pspicture}
\caption{The classification of the permutations in $S^2$ into $5$ types $A$ to $E$ \label{fig:S2-types}}
\end{center}
\end{figure}
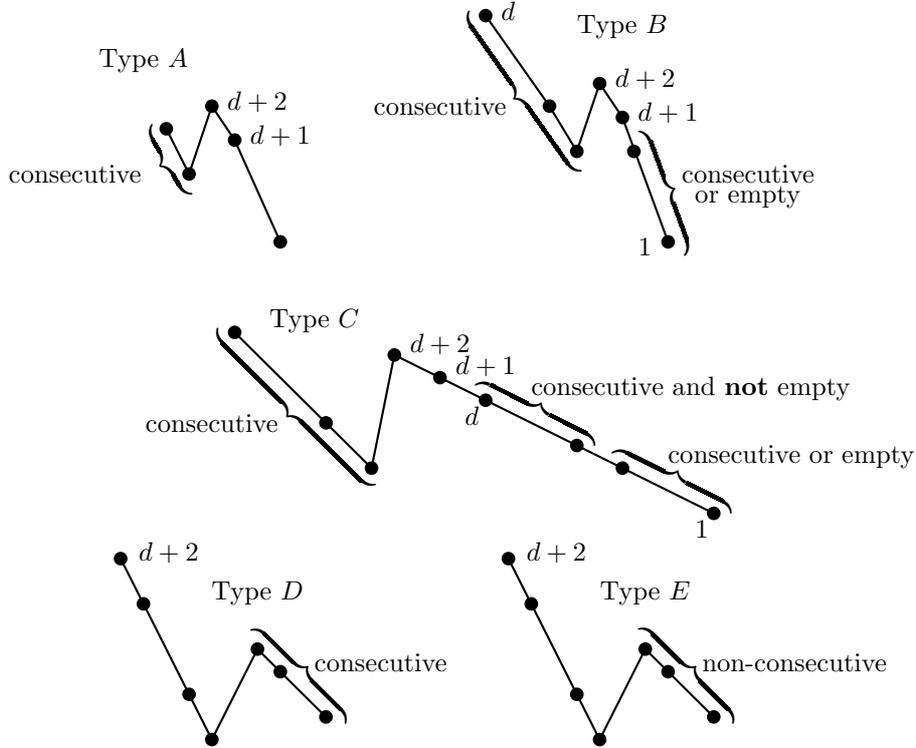

The bijection between $\mathcal{NI}$ and $S^2$ is less simple, and we will need to classify the permutations of $S^2$ by dividing them into types, from $A$ to $E$. Those types are illustrated in Figure \ref{fig:S2-types}.

The permutations $\sigma$ of type $A$ are those of $S^2$ such that (1) $d+2$ is the second element of the ascent of $\sigma$, and (2) the first sequence of descents of $\sigma$ contains only two elements, that are consecutive.

The permutations $\sigma$ of type $B$ are those of $S^2$ such that (1) $d+2$ is the second element of the ascent of $\sigma$, (2) the first sequence of descents of $\sigma$ is composed of consecutive elements, and contains at least $3$ elements, and (3) the second sequence of descents of $\sigma$ has the form $(d+2) (d+1) r$, with $r$ being either empty or a sequence of consecutive elements in decreasing order and whose smallest element is $1$.

The permutations $\sigma$ of type $C$ are those of $S^2$ such that (1) $d+2$ is the second element of the ascent of $\sigma$, (2) the first sequence of descents of $\sigma$ is made of consecutive elements, and contains at least $3$ elements, and (3) the second sequence of descents of $\sigma$ is of the form $(d+2) (d+1) r_1 r_2$ with $r_1$ being a sequence of consecutive elements in decreasing order and whose greatest element is $d$, and $r_2$ being either empty or a sequence of consecutive elements in decreasing order and whose smallest element is $1$. Notice that $r_1$ cannot be empty.

The permutations of type $D$ are those of $S^2$ such that (1) $d+2$ is the first element of $\sigma$, and (2) the elements of the second sequence of descents of $\sigma$ are consecutive.

The permutations of type $E$ are those of $S^2$ such that (1) $d+2$ is the first element of $\sigma$, and (2) the elements of the second sequence of descents of $\sigma$ are not consecutive.

Given this classification, it is now easy to prove that:
\begin{proposition}
Let $\sigma$ be a permutation of $S^2$. Then $\sigma$ is of one type exactly, among the types $A$ to $E$. \label{prop:one-type}
\end{proposition}

\begin{proof}
We distinguish two cases, according to the position of $d+2$ in $\sigma$: $d+2$ is either the first element of $\sigma$ or the second element of the ascent of $\sigma$. In the first case, it is clear that $\sigma$ is either of type $D$ or of type $E$. Let us now assume that $d+2$ is the second element of the ascent of $\sigma$. Then, because $\sigma \in S^2$, the elements of the first sequence of descents of $\sigma$ are necessarily consecutive.

Let us consider the position of $d+1$ in $\sigma$. If it is the first element of $\sigma$, and since the elements in the first sequence of descents of $\sigma$ are consecutive, then the diamond property around the ascent of $\sigma$ is not satisfied. Indeed, in such a situation, it is impossible for the rightmost element of the diamond to be greater than the lowest one. Consequently, the only possible position for $d+1$ in $\sigma$ is just after $d+2$.

If there are only two elements in the first sequence of descents of $\sigma$, then $\sigma$ is of type $A$. If there are at least three elements in the first sequence of descents of $\sigma$, then it is of type $C$ if $d+1$ is followed by $d$, of type $B$ otherwise. Because the elements in the first sequence of descents of $\sigma$ are consecutive, the reader will easily understand that the second sequence of descents of $\sigma$ is composed of consecutive element for $\sigma$ of type $B$, and splits into two sequences of consecutive elements in case $\sigma$ is of type $C$.
\end{proof}

We are now able to define the application $\Phi^2$ from $\mathcal{NI}$ to $S^2$, and to prove that it is a bijection.

\begin{figure}[ht]
\begin{center}

\psset{unit=0.3cm}
\begin{pspicture}(0,0)(40,11)
\rput[l](0.5,6){$s = \{1,2, \ldots , d+1\} \setminus \{x\}$}
\rput[l](0.5,4){$w = \{x\}$}

\rput[r](19.5,4){$\Phi^2(s) = $}

\pscircle*(22,5){0.3}
\rput[r](21.5,4.5){$x$}
\psline(22,5)(23,4)
\pscircle*(23,4){0.3}
\rput[r](22.5,3.5){$x-1$}
\psline(23,4)(25,7)
\pscircle*(25,7){0.3}
\rput[l](25.5,7.5){$d+2$}
\psline(25,7)(26,6)
\pscircle*(26,6){0.3}
\rput[l](26.5,6.5){$d+1$}
\psline(26,6)(29,3)
\pscircle*(29,3){0.3}

\rput{135}(29.2,7){$\underbrace{~ \ ~ \ \qquad \qquad \ ~}$}
\rput[l](30,7.5){all elements}
\rput[l](31.5,6){but $x$ and $x-1$}

\end{pspicture}
\caption{Definition of the bijection $\Phi^2$ for $s$ such that $|w| = 1$ \label{fig:Bijection1}}
\end{center}
\end{figure}
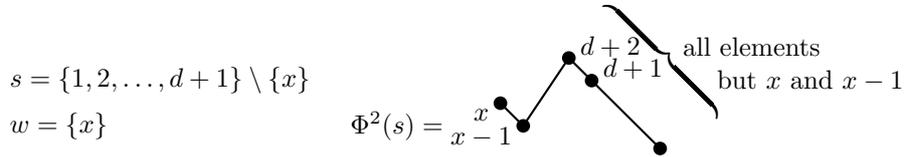

Consider a non-interval subset $s$ of $\{1,2, \ldots , d+1\}$, and call $w$ the associated set of wholes $w = \{1,2, \ldots , d+1\} \setminus s$.
\begin{enumerate}
 \item If $w$ contains only one element $x$, then necessarily $x \neq 1$ and $x \neq d+1$, or $s$ would be an interval. In this case, we set $\Phi^2(s)$ to the permutation of type $A$ with $x (x-1)$ on its first descent. This permutation obviously satisfies the diamond property (see Figure \ref{fig:Bijection1}).
\end{enumerate}

If $w$ contains at least two elements, let us denote by $n$ the cardinality of $s$ and by $m$ the cardinality of $w$ increased by $1$. Notice that $m \geq 3$ and $n \geq 2$. We will also call $w_1$ and $w_2$ the smallest and second smallest elements of $w$, and $s_n$ and $s_{n-1}$ the greatest and second greatest elements of $s$. We will associate to $s$ a permutation of $S^2$ with $m$ elements on its first sequence of descents and $n$ on its second, according to the relative order of $w_1$, $w_2$, $s_n$ and $s_{n-1}$.

Actually, there are few ways to order those $4$ elements, since they must satisfy the conditions $w_1 < w_2$, $s_{n-1} < s_n$, and $w_1 < s_n$ (or $s$ would be an interval). Namely there are five possible such orderings.

\begin{enumerate}
\addtocounter{enumi}{1}
 \item If $w_1 < w_2< s_{n-1} < s_n$ or $w_1 <s_{n-1}< w_2 < s_n $, then $\Phi^2(s)$ is the permutation of type $E$ obtained as follows: we start from $d+2$, then write the elements of $w$ in decreasing order, and finally the elements of $s$ in decreasing order. Because of the conditions satisfied by $w_1$, $w_2$, $s_n$ and $s_{n-1}$, this permutation satisfies the diamond property (see Figure \ref{fig:Bijection2}).

\begin{figure}[ht]
\begin{center}

\psset{unit=0.3cm}
\begin{pspicture}(0,0)(40,15)
\rput[l](0.5,14){$s = \{s_1, s_2, \ldots , s_n\}$ with $s_1 < s_2 < \ldots < s_n$ }
\rput[l](0.5,12){$w = \{w_1, w_2, \ldots , w_{m-1}\}$ with $w_1 < w_2 < \ldots < w_{m-1}$}
\rput[l](0.5,10){$w_1 < w_2< s_{n-1} < s_n$ or $w_1 <s_{n-1}< w_2 < s_n $}

\rput[r](6,4){$\Phi^2(s) = $}

\pscircle*(8,8){0.3}
\rput[r](7.5,7.5){$d+2$}
\psline(8,8)(10,6.5)
\pscircle*(10,6.5){0.3}
\rput[r](9.5,6){$w_{m-1}$}
\psline(10,6.5)(14,3.5)
\pscircle*(14,3.5){0.3}
\rput[r](13.5,3){$w_2$}
\psline(14,3.5)(16,2)
\pscircle*(16,2){0.3}
\rput[r](15.5,1.5){$w_1$}
\psline(16,2)(20,7)
\pscircle*(20,7){0.3}
\rput[l](20.5,7.5){$s_n$}
\psline(20,7)(22,5.5)
\pscircle*(22,5.5){0.3}
\rput[l](22.5,6){$s_{n-1}$}
\psline(22,5.5)(26,2.5)
\pscircle*(26,2.5){0.3}
\rput[l](26.5,3){$s_1$}

\rput{145}(25.2,6.3){$\underbrace{~ \ ~ \ \quad \qquad \qquad \quad \ ~}$}
\rput[l](26.5,6.5){non-consecutive}

\end{pspicture}
\caption{Definition of the bijection $\Phi^2$ for $s$ such that $w_1 < w_2< s_{n-1} < s_n$ or $w_1 <s_{n-1}< w_2 < s_n $ \label{fig:Bijection2}}
\end{center}
\end{figure}
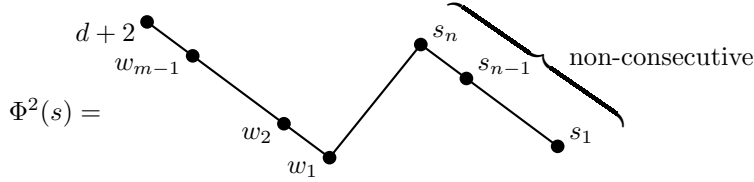

 \item If $s_{n-1} < w_1 < w_2 < s_n$, then the non-interval subset $s$ is completly determined by knowing the cardinality $n$ of $s$ and the grestest element $s_n$ of $s$. Indeed, it is necessary that $s = \{1,2, \ldots , n-1\} \uplus \{s_n\}$ to satisfy the condition $s_{n-1} < w_1 < w_2 < s_n$. In this case, we associate to $s$ a permutation of type $D$ as follows. The first element of $\Phi^2(s)$ is $d+2$, the second sequence of descents of $\Phi^2(s)$ is made of $n$ consecutive elements in decreasing order, the greatest of which is $s_n$, and the remaining elements are placed after $d+2$ in decreasing order to complete the first sequence of descents of $\Phi^2(s)$. To prove that this permutation is of type $D$, we must check that it belongs to $S^2$, that is to say that it satisfies the diamond property. It is simple to see that $s_n$ has at least $n+1$ elements smaller than itself: the remaining $n-1$ elements of $s$, $w_1$ and $w_2$. Consequently, $1$ and $2$ cannot be in the second sequence of descents of $\Phi^2(s)$. Therefore, the first sequence of descents of $\Phi^2(s)$ ends with $2 \ 1$, and this is enough to prove the diamond property (see Figure \ref{fig:Bijection3}).

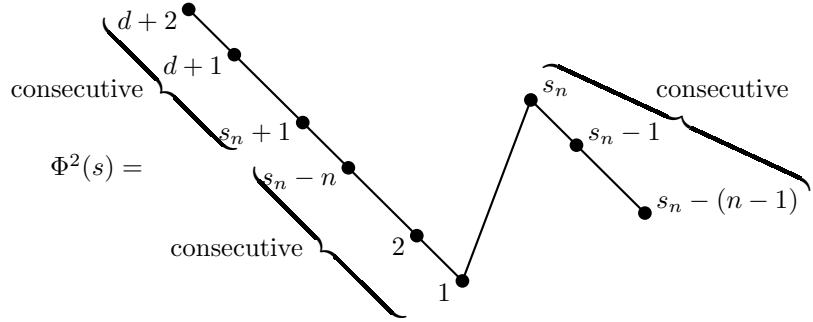
\begin{figure}[ht]
\begin{center}
\psset{unit=0.3cm}
\begin{pspicture}(0,0)(40,21)
\rput[l](0.5,20){$s_{n-1} < w_1 < w_2 < s_n$}
\rput[l](0.5,18){$s = \{1, 2, \ldots , n-1\} \uplus \{s_n\}$ }
\rput[l](0.5,16){$w = \{n,n+1 , \ldots , s_n-1\} \uplus \{s_n+1, \ldots, d+1\}$ }

\rput[r](8,7){$\Phi^2(s) = $}

\pscircle*(10,14){0.3}
\rput[r](9.5,13.5){$d+2$}
\psline(10,14)(12,12)
\pscircle*(12,12){0.3}
\rput[r](11.5,11.5){$d+1$}
\psline(12,12)(15,9)
\pscircle*(15,9){0.3}
\rput[r](14.5,8.5){$s_n+1$}
\psline(15,9)(17,7)
\pscircle*(17,7){0.3}
\rput[r](16.5,6.5){$s_n-n$}
\psline(17,7)(20,4)
\pscircle*(20,4){0.3}
\rput[r](19.5,3.5){$2$}
\psline(20,4)(22,2)
\pscircle*(22,2){0.3}
\rput[r](21.5,1.5){$1$}
\psline(22,2)(25,10)
\pscircle*(25,10){0.3}
\rput[l](25.5,10.5){$s_n$}
\psline(25,10)(27,8)
\pscircle*(27,8){0.3}
\rput[l](27.5,8.5){$s_n-1$}
\psline(27,8)(30,5)
\pscircle*(30,5){0.3}
\rput[l](30.5,5.5){$s_n-(n-1)$}

\rput{-45}(9,10.5){$\underbrace{~ \ ~ \ \qquad \qquad \quad \ ~}$}
\rput[r](8,10.5){consecutive}
\rput{-45}(16,3.5){$\underbrace{~ \ ~ \ \qquad \qquad \qquad \ ~}$}
\rput[r](15,3.5){consecutive}

\rput{155}(31.5,9){$\underbrace{~ \ ~ \ \quad \qquad \qquad \qquad \qquad \ ~}$}
\rput[l](30.5,10.5){consecutive}

\end{pspicture}
\caption{Definition of the bijection $\Phi^2$ for $s$ such that $s_{n-1} < w_1 < w_2 < s_n$ \label{fig:Bijection3}}
\end{center}
\end{figure}

 \item If $w_1 < s_{n-1} < s_n < w_2$, the elements of $\{1,2, \ldots , d+1\}$ are partitioned into $s \uplus w$ in the following way : $s = \{1,\ldots, w_1 -1\} \uplus \{w_1 +1, \ldots , n+1\}$ and $w = \{w_1\} \uplus \{w_2 = n+2, \ldots, d+1\}$. The non-interval $s$ is then completly determined by knowing the cardinality $n$ of $s$ and the number $p = n+1-w_1$ of elements of $s$ between $w_1$ and $w_2$. Let us notice that $p \geq 2$ (since $s_{n-1}$ and $s_n$ are between $w_1$ and $w_2$) and $p \leq n-1$ ($p=n$ would imply that $s$ is an interval). In this case, we associate to $s$ the permutation $\Phi^2(s)$ of type $C$ as follows. The second sequence of descents of $\Phi^2(s)$ splits into two parts (the second one possibly empty). The first part contains $p+1$ elements (we can check that $3 \leq p+1 \leq n$) that are consecutive, and whose greatest element is $d+2$, of course written in decreasing order. The second part is composed of $n-p-1$ consecutive elements in decreasing order, with $1$ as minimal element. This construction leaves $m$ consecutive elements unused so far: written in decreasing order, they will constitute the first sequence of descents of $\Phi^2(s)$. Now, it is easy to prove the diamond property, since the second sequence of descents of $\Phi^2(s)$ necessarily starts with $(d+2)(d+1)$. This remark completes the proof that the permutation $\Phi^2(s)$ we just defined is in $S^2$, and of type $C$ (see Figure \ref{fig:Bijection4}).

\begin{figure}[ht]
\begin{center}
\psset{unit=0.3cm}
\begin{pspicture}(0,0)(40,21)
\rput[l](0.5,20){$ w_1 < s_{n-1} < s_n < w_2 $}
\rput[l](0.5,18){$s = \{1, \ldots , w_1-1\} \uplus \{w_1+1, \ldots , n+1\}$ }
\rput[l](0.5,16){$w = \{w_1\} \uplus \{n+2, \ldots, d+1\}$ }
\rput[l](0.5,14){We set $p = n+1-w_1$}

\rput[r](5,10){$\Phi^2(s) = $}

\pscircle*(10,9){0.3}
\rput[r](9.5,8.5){$d+1-p$}
\psline(10,9)(13,6)
\pscircle*(13,6){0.3}
\rput[r](12.5,5.5){$n-p+1$}
\psline(13,6)(15,4)
\pscircle*(15,4){0.3}
\rput[r](14.5,3.5){$n-p$}
\psline(15,4)(18,12)
\pscircle*(18,12){0.3}
\rput[l](18.5,12.5){$d+2$}
\psline(18,12)(20,10)
\pscircle*(20,10){0.3}
\rput[l](20.5,10.5){$d+1$}
\psline(20,10)(23,7)
\pscircle*(23,7){0.3}
\rput[l](23.5,7.5){$d+2-p$}
\psline(23,7)(25,5)
\pscircle*(25,5){0.3}
\rput[l](25.5,5.5){$n-p-1$}
\psline(25,5)(28,2)
\pscircle*(28,2){0.3}
\rput[l](28.5,2.5){$1$}

\rput{135}(24,11.5){$\underbrace{~ \ ~ \ \qquad \qquad \qquad \ ~}$}
\rput[l](25,11.5){$p+1$ consecutive elements}
\rput{135}(31.5,5){$\underbrace{~ \ ~ \ \qquad \quad \ ~}$}
\rput[l](31,6){$n-p-1$}
\rput[l](32.2,5){consecutive}
\rput[l](33.2,4){elements}

\rput{-45}(7.2,5.5){$\underbrace{~ \ ~ \ \quad \qquad \qquad \ ~}$}
\rput[r](5.5,6.5){$m$}
\rput[r](6.5,5.5){consecutive}
\rput[r](7.5,4.5){elements}

\end{pspicture}
\caption{Definition of the bijection $\Phi^2$ for $s$ such that $w_1 < s_{n-1} < s_n < w_2$ \label{fig:Bijection4}}
\end{center}
\end{figure}
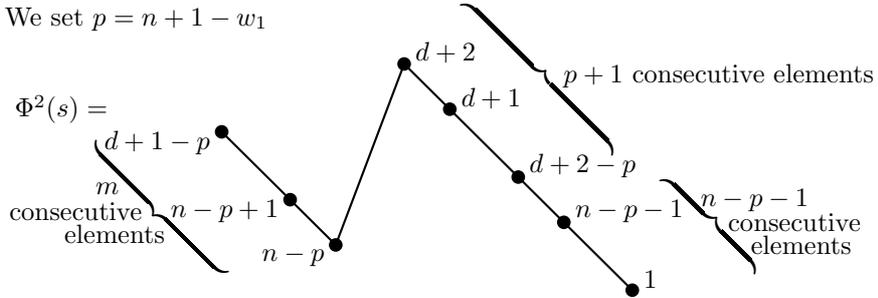

 \item The last possible relative order of $w_1$, $w_2$, $s_n$ and $s_{n-1}$ is $s_{n-1} < w_1 < s_n < w_2$. This case is particularly simple since the cardinality $n$ of $s$ determines $s$ completly. Indeed, it is necessary that $s=\{1, \ldots, n-1\} \uplus \{n+1\}$ to satisfy the conditions $s_{n-1} < w_1 < s_n < w_2$. The permutation $\Phi^2(s)$ is of type $B$, with the $n$ elements on the second sequence of descents starting with $(d+2)(d+1)$ and then either nothing or consecutive numbers in decreasing order and ending with $1$. This leaves $m$ consecutive numbers, with greatest element $d$, to fill in the first sequence of descents of $\Phi^2(s)$. Because the second sequence of descents starts with $(d+2)(d+1)$, $\Phi^2(s)$ clearly satisfies the diamond property, justifying that $\Phi^2(s)$ is a permutation of $S^2$ and of type $B$ (see Figure \ref{fig:Bijection5}).

\begin{figure}[ht]
\psset{unit=0.3cm}
\begin{pspicture}(0,4)(40,21)
\rput[l](0.5,20){$s_{n-1} < w_1 < s_n < w_2$}
\rput[l](0.5,18){$s = \{1, \ldots , n-1\} \uplus \{n+1\}$ }
\rput[l](0.5,16){$w = \{n\} \uplus \{n+2, \ldots, d+1\}$ }

\rput[r](5,9){$\Phi^2(s) = $}

\pscircle*(12,12){0.3}
\rput[r](11.5,11.5){$d$}
\psline(12,12)(14,10)
\pscircle*(14,10){0.3}
\rput[r](13.5,9.5){$d-1$}
\psline(14,10)(17,7)
\pscircle*(17,7){0.3}
\rput[r](16.5,6.5){$n$}
\psline(17,7)(19,5)
\pscircle*(19,5){0.3}
\rput[r](18.5,4.5){$n-1$}
\psline(19,5)(22,14)
\pscircle*(22,14){0.3}
\rput[l](22.5,14.5){$d+2$}
\psline(22,14)(24,12)
\pscircle*(24,12){0.3}
\rput[l](24.5,12.5){$d+1$}
\psline(24,12)(29,7)
\pscircle*(26,10){0.3}
\rput[l](26.5,10.5){$n-2$}
\psline(26,10)(31,5)
\pscircle*(31,5){0.3}
\rput[l](31.5,5.5){$1$}

\rput{135}(31.5,9){$\underbrace{~ \ ~ \ \qquad \qquad \quad \ ~}$}
\rput[l](33,10){$n-2$}
\rput[l](32.2,9){consecutive}
\rput[l](33.2,8){elements}

\rput{-45}(12.2,7.5){$\underbrace{~ \ ~ \ \quad \qquad \qquad \qquad \ ~}$}
\rput[r](10.5,8.5){$m$}
\rput[r](11.5,7.5){consecutive}
\rput[r](12.5,6.5){elements}

\end{pspicture}
\caption{Definition of the bijection $\Phi^2$ for $s$ such that $s_{n-1} < w_1 < s_n < w_2$ \label{fig:Bijection5}}
\end{figure}
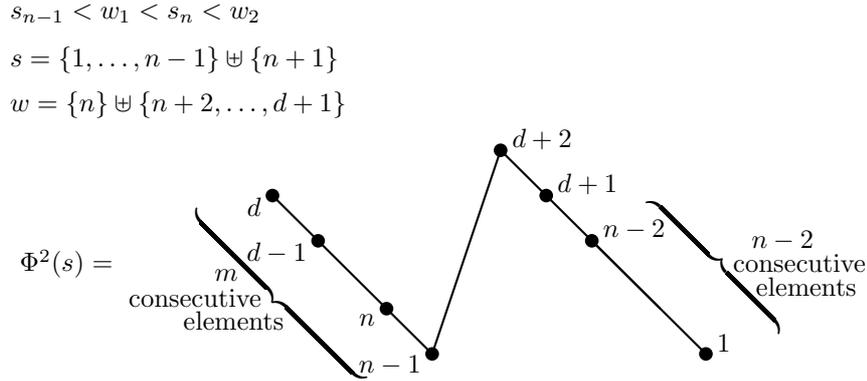
\end{enumerate}

These different cases to define $\Phi^2(s)$ are exemplified in Figure \ref{tab:phi2}.

\begin{figure}[ht]
\begin{center}
\psset{unit=0.19cm}
\begin{tabular}{|p{3.3cm}|p{2.5cm}|c|c|}
\hline
Case for $s$ & Example of $s$ & $\Phi^2(s)$ & Type\\
\hline
(1) with & $s=\{1,2,4,5,6\}$ & \multirow{3}{*}{\begin{pspicture}(-1,-1)(10,7)
\pscircle*(1,3){0.4}
\rput(0,2){$3$}
\psline(1,3)(2,2)
\pscircle*(2,2){0.4}
\rput(1,1){$2$}
\psline(2,2)(4,5)
\pscircle*(4,5){0.4}
\rput(5,6){$7$}
\psline(4,5)(5,4)
\pscircle*(5,4){0.4}
\rput(6,5){$6$}
\psline(5,4)(6,3)
\pscircle*(6,3){0.4}
\rput(7,4){$5$}
\psline(6,3)(7,2)
\pscircle*(7,2){0.4}
\rput(8,3){$4$}
\psline(7,2)(8,1)
\pscircle*(8,1){0.4}
\rput(9,2){$1$}
\end{pspicture}} & \multirow{3}{*}{$A$}  \\
$w = \{x\}$ & $w=\{3\}$ & & \\
 & $x=3$ & & \\
\hline
(2)  with & $s=\{1,5,6\}$ & \multirow{3}{*}{\begin{pspicture}(-1,-1)(10,6)
\pscircle*(1,4){0.4}
\rput(0,3){$7$}
\psline(1,4)(2,3)
\pscircle*(2,3){0.4}
\rput(1,2){$4$}
\psline(2,3)(3,2)
\pscircle*(3,2){0.4}
\rput(2,1){$3$}
\psline(3,2)(4,1)
\pscircle*(4,1){0.4}
\rput(3,0){$2$}
\psline(4,1)(6,3.5)
\pscircle*(6,3.5){0.4}
\rput(7,4.5){$6$}
\psline(6,3.5)(7,2.5)
\pscircle*(7,2.5){0.4}
\rput(8,3.5){$5$}
\psline(7,2.5)(8,1.5)
\pscircle*(8,1.5){0.4}
\rput(9,2.5){$1$}
\end{pspicture}} & \multirow{3}{*}{$E$}  \\
$w_1 < w_2 < s_{n-1} < s_n$ & $w=\{2,3,4\}$ & & \\
 & & & \\
\hline
(2)  with & $s=\{1,3,4,6\}$ & \multirow{3}{*}{\begin{pspicture}(-1,-1)(10,6)
\pscircle*(1,3.5){0.4}
\rput(0,2.5){$7$}
\psline(1,3.5)(2,2.5)
\pscircle*(2,2.5){0.4}
\rput(1,1.5){$5$}
\psline(2,2.5)(3,1.5)
\pscircle*(3,1.5){0.4}
\rput(2,0.5){$2$}
\psline(3,1.5)(5,4)
\pscircle*(5,4){0.4}
\rput(6,5){$6$}
\psline(5,4)(6,3)
\pscircle*(6,3){0.4}
\rput(7,4){$4$}
\psline(6,3)(7,2)
\pscircle*(7,2){0.4}
\rput(8,3){$3$}
\psline(7,2)(8,1)
\pscircle*(8,1){0.4}
\rput(9,2){$1$}
\end{pspicture}} &  \multirow{3}{*}{$E$} \\
$w_1 < s_{n-1} < w_2 < s_n$ & $w=\{2,5\}$ & & \\
 & & & \\
\hline
(3) with & $s=\{1,2,5\}$ & \multirow{3}{*}{\begin{pspicture}(-1,-1)(10,6)
\pscircle*(1,4){0.4}
\rput(0,3){$7$}
\psline(1,4)(2,3)
\pscircle*(2,3){0.4}
\rput(1,2){$6$}
\psline(2,3)(3,2)
\pscircle*(3,2){0.4}
\rput(2,1){$2$}
\psline(3,2)(4,1)
\pscircle*(4,1){0.4}
\rput(3,0){$1$}
\psline(4,1)(6,3.5)
\pscircle*(6,3.5){0.4}
\rput(7,4.5){$5$}
\psline(6,3.5)(7,2.5)
\pscircle*(7,2.5){0.4}
\rput(8,3.5){$4$}
\psline(7,2.5)(8,1.5)
\pscircle*(8,1.5){0.4}
\rput(9,2.5){$3$}
\end{pspicture}} &  \multirow{3}{*}{$D$} \\
 $s_{n-1} < w_1 < w_2 < s_n$ & $w=\{3,4,6\}$ & & \\
 & $|s| = 3$, $s_n=5$ & & \\
\hline
(4) with  & $s=\{1,3,4,5\}$ & \multirow{3}{*}{\begin{pspicture}(-1,-1)(10,6)
\pscircle*(1,3.5){0.4}
\rput(0,2.5){$3$}
\psline(1,3.5)(2,2.5)
\pscircle*(2,2.5){0.4}
\rput(1,1.5){$2$}
\psline(2,2.5)(3,1.5)
\pscircle*(3,1.5){0.4}
\rput(2,0.5){$1$}
\psline(3,1.5)(5,4)
\pscircle*(5,4){0.4}
\rput(6,5){$7$}
\psline(5,4)(6,3)
\pscircle*(6,3){0.4}
\rput(7,4){$6$}
\psline(6,3)(7,2)
\pscircle*(7,2){0.4}
\rput(8,3){$5$}
\psline(7,2)(8,1)
\pscircle*(8,1){0.4}
\rput(9,2){$4$}
\end{pspicture}} &  \multirow{3}{*}{$C$} \\
$w_1 < s_{n-1} < s_n < w_2$ & $w=\{2,6\}$ & & \\
 & $|s|=4$, $p=3$ & & \\
\hline
(4) with  & $s=\{1,2,4,5\}$ & \multirow{3}{*}{\begin{pspicture}(-1,-1)(10,6)
\pscircle*(1,3.5){0.4}
\rput(0,2.5){$4$}
\psline(1,3.5)(2,2.5)
\pscircle*(2,2.5){0.4}
\rput(1,1.5){$3$}
\psline(2,2.5)(3,1.5)
\pscircle*(3,1.5){0.4}
\rput(2,0.5){$2$}
\psline(3,1.5)(5,4)
\pscircle*(5,4){0.4}
\rput(6,5){$7$}
\psline(5,4)(6,3)
\pscircle*(6,3){0.4}
\rput(7,4){$6$}
\psline(6,3)(7,2)
\pscircle*(7,2){0.4}
\rput(8,3){$5$}
\psline(7,2)(8,1)
\pscircle*(8,1){0.4}
\rput(9,2){$1$}
\end{pspicture}} &  \multirow{3}{*}{$C$} \\
$w_1 < s_{n-1} < s_n < w_2$ & $w=\{3,6\}$ & & \\
 & $|s|=4$, $p=2$ & & \\
\hline
(5) with  & $s=\{1,2,3,5\}$ & \multirow{3}{*}{\begin{pspicture}(-1,-1)(10,6)
\pscircle*(1,3.5){0.4}
\rput(0,2.5){$5$}
\psline(1,3.5)(2,2.5)
\pscircle*(2,2.5){0.4}
\rput(1,1.5){$4$}
\psline(2,2.5)(3,1.5)
\pscircle*(3,1.5){0.4}
\rput(2,0.5){$3$}
\psline(3,1.5)(5,4)
\pscircle*(5,4){0.4}
\rput(6,5){$7$}
\psline(5,4)(6,3)
\pscircle*(6,3){0.4}
\rput(7,4){$6$}
\psline(6,3)(7,2)
\pscircle*(7,2){0.4}
\rput(8,3){$2$}
\psline(7,2)(8,1)
\pscircle*(8,1){0.4}
\rput(9,2){$1$}
\end{pspicture}} &  \multirow{3}{*}{$B$} \\
$s_{n-1} < w_1 < s_n < w_2$ & $w=\{4,6\}$ & & \\
 & $|s|=4$ & & \\
\hline
\end{tabular}
\caption{Definition of $\Phi^2(s)$ for some non-interval subsets $s$ of $\{1,2,\ldots,6\}$ ($d=5$), illustrating all the possible cases in the construction of $\Phi^2$ \label{tab:phi2}}
\end{center}
\end{figure}

This ends the definition of the application $\Phi^2 : \mathcal{NI} \rightarrow S^2$. Moreover, we have:
\begin{proposition}
The application $\Phi^2$ defines a bijection between $\mathcal{NI}$ and $S^2$. \label{prop:phi2}
\end{proposition}

\begin{proof}
The inverse application of $\Phi^2$, from $S^2$ to $\mathcal{NI}$, can easily be defined from the previous paragraphs, distinguishing cases according to the type (from $A$ to $E$) of a permutation of $S^2$. The details are left to the reader.
\end{proof}

Putting things all together, we have a partition of the set of minimal permutations with $d$ descents and of size $d+2$ into $S^1 \uplus S^2$, and two bijections $\Phi^1$ (resp. $\Phi^2$) between $S^1$ (resp. $S^2$) and $\mathcal{NI}$. Combining this with the enumeration of non-interval subsets of $\{1,2, \ldots, d+1\}$ obtained in Proposition \ref{prop:non-interval}, we get another proof of Theorem \ref{thm:d+2}, by a bijective approach.

\section{Conclusion and open problems}
\label{section:conclusion}

The goal pursued in this paper is the analysis (characterization, enumeration, $\ldots$) of the permutations that are minimal for the property of having $d$ descents, minimal being intended in the sense of the pattern-involvement relation. For $d=2^p$, those permutations arise from the whole genome duplication - random loss model, defined in computational biology, where they appear as the excluded patterns defining the pattern-avoiding classes of permutations obtained in at most $p$ steps in this model.

We first provided a local characterization of the minimal permutations with $d$ descents, focusing only on the elements of the permutation surrounding its ascents. This characterization is easy to check: indeed, it provides a linear-time procedure for deciding whether a permutation is minimal with $d$ descents or not.

The second step of our study was more about enumerating these permutations. We proved that a minimal permutation with $d$ descents has size at least $d+1$ and at most $2d$. We could not find the enumeration of all minimal permutations with $d$ descents, but we were able to enumerate such permutations of size $d+1$, $d+2$ and $2d$. More precisely, there is only one of size $d+1$ (which is the reversed identity), there are $2^{d+2} - (d+1)(d+2) -2$ minimal permutations with $d$ descents of size $d+2$, and those of size $2d$ are enumerated by the Catalan numbers.

The enumeration of the minimal permutations with $d$ descents and of size $n \in [(d+3)..(2d-1)]$ remains an open question. For $n=d+3$, we computed the first few terms of the enumerating sequence, and it seems not to appear in the Online Encyclopedia of Integer Sequences \cite{njas}. Notice however that the analytical technique used to enumerate the minimal permutations with $d$ descents of size $d+2$ could theoretically be applied to any other size $n \in [(d+3)..(2d-1)]$, but there would be many more cases to consider. Indeed, only for $n=d+3$, there are more than eighty of them, instead of the two cases for $n=d+2$. This combinatorial complexity suggests that to solve this enumerating problem, either other techniques or an automated examination of the numerous cases are needed.

\nocite{*}
\bibliographystyle{plain}
\bibliography{biblio}

\end{document}